\documentclass[11pt,a4paper]{article}
\usepackage{amsmath, amssymb, theorem, latexsym,graphicx}
\usepackage{fullpage}
\usepackage{color}
\usepackage[german,english]{babel}
\newtheorem{theorem}{Theorem}[section]
\newtheorem{lemma}[theorem]{Lemma}
\newtheorem{proposition}[theorem]{Proposition}
\newtheorem{definition}[theorem]{Definition}

\newtheorem{property}[theorem]{Property}
        \newtheorem{properties}[theorem]{Properties}

\numberwithin{equation}{section}
\numberwithin{figure}{section}

\allowdisplaybreaks

\usepackage[pdftex]{hyperref}

\def\Dg {{\cal D}} 
\def\Fg {{\cal F}}
\def\Lg {{\cal L}} 
\def\Qg {{\cal Q}} 
\def\Sg {{\cal S}} 
\def\Tg {{\cal T}} 
\def\Wg {{\cal W}} 

\def \nzb {{\mathbb N\setminus\{0\}}}

\def \rz {{\mathbb R}}

\def \T {{\mathbb T}}

\def\la {{\lambda}}

\newcommand {\ar}{\rightarrow}

\newcommand{\noib}{\noindent $\bullet $~}
\newcommand{\ie}{i.e., }
\newcommand{\eg}{e.g.}

\newcommand{\resp}{\emph{resp. }}

\newcommand{\pid}{\frac{\pi}{2}}
\newcommand{\piq}{\frac{\pi}{4}}
\newcommand{\pitq}{\frac{3\pi}{4}}

\title{Dirichlet eigenfunctions of the square membrane:\\
Courant's property, and A.~Stern's  and  {\AA}.~Pleijel's analyses}
\author{P. B\'erard \\ Institut Fourier, Universit\'{e} de Grenoble and CNRS, B.P.74,
\\ F 38402 Saint Martin d'H\`{e}res Cedex, France.\\
pierrehberard@gmail.com\\
and \\
B. Helffer \\
Laboratoire de Math\'ematiques, Univ. Paris-Sud and CNRS,\\
F 91405 Orsay Cedex, France, and\\
Laboratoire Jean Leray, Universit\'{e} de Nantes and CNRS\\
F 44322 Nantes Cedex 3, France.\\
Bernard.Helffer@math.u-psud.fr}

\date{\small{March 11, 2015}}

\begin{document}
\maketitle

\begin{quote}
To appear in Springer Proceedings in Mathematics \& Statistics  (2015), MIMS-GGTM conference in memory of M.~S.~Baouendi. Ali Baklouti, Aziz El Kacimi, Sadok Kallel, and Nordine Mir Editors.
\end{quote}

\bigskip
{\large{To the memory of M. Salah Baouendi}}
\bigskip

\begin{abstract}
In this paper, we revisit Courant's nodal domain theorem for the
Dirichlet eigenfunctions of a square membrane, and the analyses of
A.~Stern and {\AA}.~Pleijel.
\end{abstract}



\section{Introduction}\label{S-HI}

Courant's celebrated nodal domain theorem \cite{Cou} says that
the number of nodal domains of an eigenfunction associated with a
$k$-th eigenvalue of the Dirichlet Laplacian is less than or equal
to $k$. Here, the eigenvalues are chosen to be positive, and listed
in increasing order. It follows from a theorem of Pleijel \cite{Pl}
that equality in Courant's theorem only occurs for finitely many
values of $k$. In this case, we speak of the Courant sharp
situation. We refer to \cite{HHO, HHOT} for the connection of this
property with the question of minimal spectral partitions. \medskip

In the case of the square, it is immediate that the first, second
and fourth eigenvalues are Courant sharp. In the first part of this
note, Sections~\ref{S-HPA} and \ref{S-H3}, we provide some missing
arguments in Pleijel's paper leading to the conclusion that there
are no other cases.\medskip

In the second part of this paper, we discuss some results of Antonie
Stern. She was a PhD student of R. Courant, and defended her PhD in
1924, see \cite{St, St1}, and \cite[p. 180]{Vo}. In her thesis, she
in particular provides an infinite sequence of Dirichlet
eigenfunctions for the square, as well as an infinite sequence of
spherical harmonics on the $2$-sphere, which have exactly two nodal
domains. In this paper, we focus on her results concerning the
square, and refer to \cite{BeHe2} for an analysis of the spherical
case. In Section~\ref{HStern}, we analyze Stern's argument, leading
to the conclusion that her proofs are not quite complete. In
Section~\ref{S-1R}, we provide a detailed proof of Stern's main
result for the square, Theorem~\ref{stern-T}.

\medskip

The authors would like to thank Virginie Bonnaillie-No\"{e}l for her
pictures of nodal domains \cite{vbn}, and Annette Vogt for her
biographical information on A.~Stern.  The authors are indebted to
D.~Jakobson for pointing out the unpublished report \cite{GSP}, and
to M.~Persson-Sundqvist for useful remarks. The authors thank
the anonymous referee for his comments and careful reading. The
second author would like to thank T.~Hoffmann-Ostenhof for
motivating discussions.


\section{Pleijel's analysis}\label{S-HPA}

Consider the rectangle $\mathcal R(a,b)=]0,a\pi[\times ]0,b\pi[$.
The  Dirichlet eigenvalues for $-\Delta$ are given by
$$\hat \la_{m,n}=(\frac{m^2}{a^2} +\frac{n^2}{b^2})\,, ~ m,n \ge 1,$$
with a corresponding basis of eigenfunctions given by
$$
\phi_{m,n} (x,y)= \sin \frac{mx}{a} \,\sin \frac{ny}{b}\,.
$$

It is easy to determine the Courant sharp eigenvalues when $b^2/a^2$
is irrational (see for example  \cite{HHOT}). The rational case is
more difficult.  Let us analyze the zero set of the Dirichlet
eigenfunctions for the square. If we normalize the square as
$]0,\pi[\times ]0,\pi[$, we have,
$$
\phi_{m,n} (x,y) = \phi_m(x)\phi_n (y)\,,\, \mbox{ with }
\phi_m(t)=\sin (m t)\,.
$$
Due to multiplicities, we have (at least) to consider the family of
eigenfunctions,
$$(x,y) \mapsto \Phi_{m,n}(x,y,\theta):= \cos \theta\, \phi_{m,n}(x,y)
+ \sin \theta\, \phi_{n,m}(x,y)\,,
$$
with $m, n \ge 1$, and $\theta \in [0,\pi[$.

In \cite{Pl}, Pleijel claims that the Dirichlet eigenvalue
$\lambda_k$ of the square is Courant sharp if and only if $k=1,2,4$.
The key point in his proof is to exclude the eigenvalues $\lambda_5,
\lambda_7$ and $\lambda_9$ which correspond respectively to the
pairs $(m,n)= (1,3)$,  $(m,n)= (2,3)$ and $(m,n)=(1,4)$.

\medskip



\begin{figure}[!ht]
\begin{center}
\includegraphics[width=12cm]{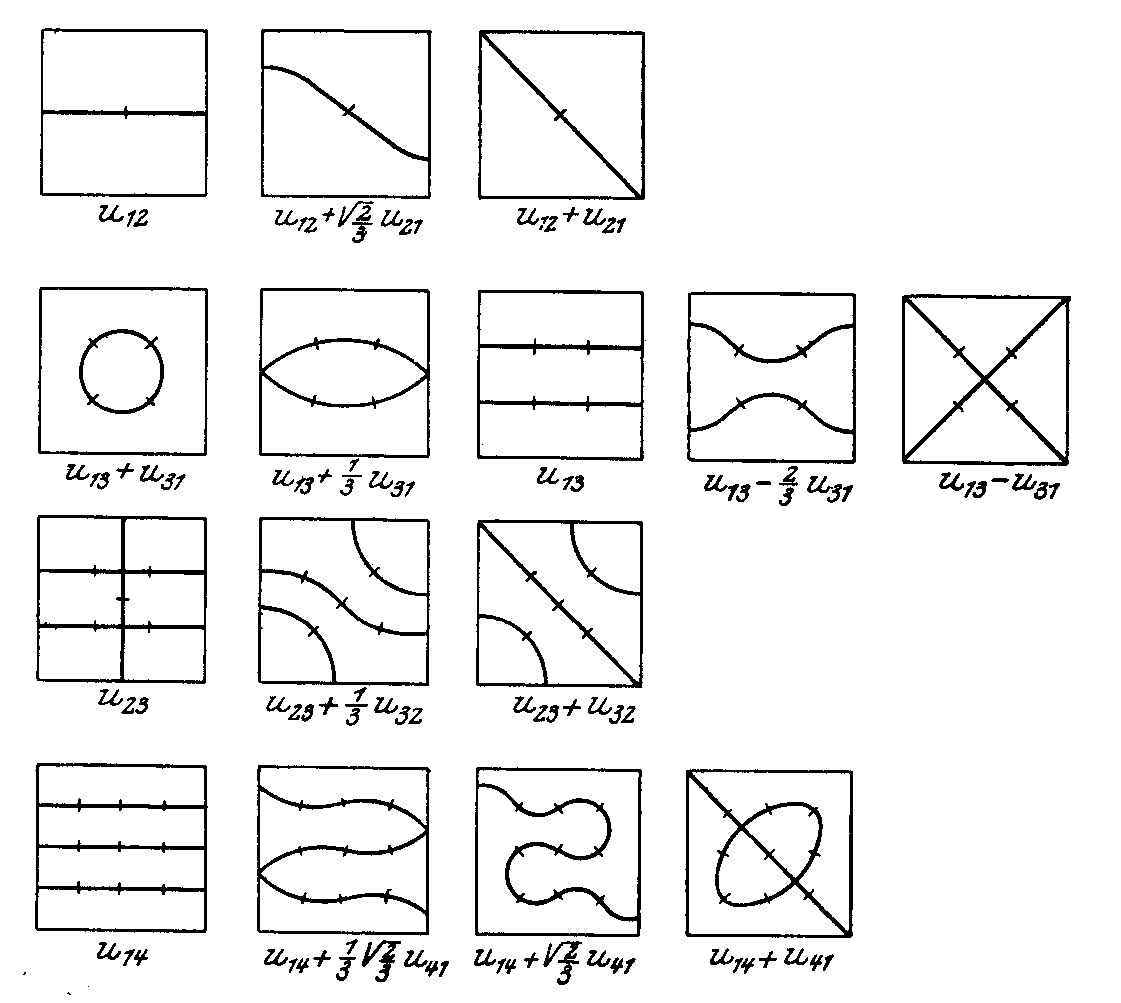}
\caption{Nodal sets for the Dirichlet eigenvalues $\lambda_2$,
$\lambda_5$, $\lambda_7$ and $\lambda_9$
(reproduced from \cite{CHd})}\label{F-pch}
\end{center}
\end{figure}

Let us briefly recall Pleijel's argument. Let $N(\lambda) :=
\#\left\{n ~|~ \lambda_n < \lambda \right\}$ be the counting
function. Using a covering of $\rz^2$ by the squares
$]k,k+1[\times]\ell,\ell+1[$, he first establishes the
estimate\footnote{There is an unimportant sign error in \cite{Pl}.}
\begin{equation}
N(\lambda) > \frac \pi 4 \lambda - 2 \sqrt{\lambda} - 1\,.
\end{equation}

 If $\lambda_n$ is Courant sharp, then $\lambda_{n-1} <
\lambda_n$, hence $N(\lambda_n)= n-1$, and
\begin{equation}\label{i1}
n >\frac \pi 4 \lambda_n - 2 \sqrt{\lambda_n}\,.
\end{equation}
On the other hand, if $\lambda_n$ is Courant sharp, the Faber-Krahn
inequality  \cite{Ban, EM-RFK} gives the necessary condition
\begin{equation*}
\frac{\lambda_n}{n} \geq \frac{{j_{0,1}}^2}{\pi}
\end{equation*}
or
\begin{equation}\label{i2}
\frac{n}{\lambda_n} \leq \pi {j_{0,1}}^{-2} <  0.54323\,.
\end{equation}
Recall that $j_{0,1}$ is the smallest positive zero of the
Bessel function of order $0$, and that $ \pi {j_{0,1}}^2$ is
the ground state energy of the disk of area $1$.

Combining \eqref{i1} and \eqref{i2}, leads to the inequality

\begin{equation}\label{i3}
\lambda_n \le 68\,.
\end{equation}

After re-ordering the values $m^2+n^2$, we get the following
spectral sequence for $\lambda_n \le 73$,
\begin{equation}\label{liste}
\begin{array}{cccc}
\lambda_1 = 2 \,,& \lambda_2 = \lambda_3 = 5\,,& \lambda_4 = 8 \,,& \lambda_5 = \lambda_6 = 10
\,,\\
\lambda_{7} = \lambda_{8} = 13\,,& \lambda_{9} = \lambda_{10} =
17\,,& \lambda_{11} = 18 \,,& \lambda_{12} = \lambda_{13} = 20\,,\\
\lambda_{14} = \lambda_{15} = 25\,,& \lambda_{16} = \lambda_{17} =
26\,,& \lambda_{18} = \lambda_{19}
= 29\,,& \lambda_{20} = 32\,,\\
\lambda_{21} = \lambda_{22} = 34\,,& \lambda_{23} = \lambda_{24} =
37 \,,& \lambda_{25} = \lambda_{26}
= 40\,,& \lambda_{27} = \lambda_{28} = 41\,,\\
\lambda_{29} = \lambda_{30} = 45 \,,& \lambda_{31} = \lambda_{32} =
\lambda_{33} = 50\,,& \lambda_{34} = \lambda_{35}
= 52\,,& \lambda_{36} = \lambda_{37} = 53\,,\\
\lambda_{38} = \lambda_{39} = 58 \,,& \lambda_{40} = \lambda_{41} =
61 \,,& \lambda_{42} = \lambda_{43} = 65\,,& \lambda_{44} = \lambda_{45} = 65\,,\\
\lambda_{46} = \lambda_{47} = 68\,,& \lambda_{48} = 72\,,&
\lambda_{49} = \lambda_{50} = 73\,,& \cdots \,.\\
\end{array}
\end{equation}

It remains to analyze, among the eigenvalues which are less than or
equal to $ 68$, those which satisfy \eqref{i2}, and hence which
can be Courant sharp. Computing the quotients $\frac{n}{\lambda_n} $
in the list \eqref{liste}, leaves us with the eigenvalues
$\lambda_5$, $\lambda_7$ and $\lambda_9$. For these last three
cases, Pleijel refers to pictures in Courant-Hilbert \cite{CH}, \S
V.5.3 p. 302, actually reproduced from \cite{Poc}, \S II.B.6, p.~80,
see Figure~\ref{F-pch} in which $u_{mn}(x,y)$ stands for $\sin(mx) \sin(ny)$.  Although the details
are not provided in these textbooks, the choice of the parameter
values in the pictures suggests that some theoretical analysis is
involved. It is however not clear whether the displayed nodal
patterns represent all possible shapes, up to deformation or
symmetries. The difficulty being that we have to analyze the nodal
sets of eigenfunctions living in two-dimensional eigenspaces.
Clearly, $\phi_{m,n} $ has $mn$ nodal components. This corresponds
to the ``product'' situation with $\theta=0$ or $\theta = \frac \pi
2$. The figures illustrate the fact that the number of nodal domains
for a linear combination of two given independent eigenfunctions can
be smaller or larger than the number of nodal domains of the given
eigenfunctions. For these reasons, Pleijel's argument does not
appear fully convincing. In Section~\ref{S-H3}, we give a detailed
proof that eigenvalues $\lambda_5$, $\lambda_7$ and $\lambda_9$ are
not Courant sharp.\medskip

 \textbf{Remark}. As pointed out to us by
M.~Persson-Sundqvist, the last two cases can be easily dealt with
using the following adaptation of an observation due to J.~Leydold
\cite{Ley}. For any $(m,n)$,
$$
\Phi_{m,n}(\pi-x,\pi-y,\theta) =
(-1)^{m+n}\Phi_{m,n} (x,y,\theta)\,.
$$
When $(m+n)$ is odd, the corresponding eigenfunction function is odd
under the symmetry with respect to the center of the square,
and hence must have an even number of nodal domains. Our results are
actually stronger, and describe the variation of the nodal sets.

\section{The three remaining cases of Pleijel}\label{S-H3}

Behind all the computations, we have the property that, for $x\in
]0,\pi[$,
\begin{equation}
\sin mx = \sqrt{1-u^2} \, U_{m-1} (u)\,,
\end{equation}
where $U_{m-1}$ is the Chebyshev polynomial of second type and
$u=\cos x$, see \cite{MOS}.

\subsection{First case : eigenvalue $\lambda_5$, or $(m,n)=(1,3)$.}

We look at the zeroes of $\Phi_{1,3}(x,y,\theta)$, see
Figure~\ref{F-pch}, 2nd row. Since
$\Phi_{1,3}(x,y,\pid-\theta) = \Phi_{1,3}(y,x,\theta)$, we can
reduce the analysis of the nodal patterns to $\theta \in
[\piq,\pitq]$. Let,
\begin{equation}\label{eq0}
\cos x = u\,,\, \cos y=v\,.
\end{equation}
This is a $C^\infty$ change of variables from the square
$]0,\pi[\times ]0,\pi[$ onto $]-1,+1[\times ]-1,+1[$.  In these
coordinates, the zero set  of  $\Phi_{1,3}(x,y,\theta)$ inside the
square is given by
\begin{equation}\label{eq1}
\cos \theta \, ( 4v^2 - 1) + \sin \theta\,  (4u^2 - 1) =0\,.
\end{equation}

To completely determine the nodal set, we have to take the closure
in $[-1,1]\times [-1,1]$ of the zero set \eqref{eq1}. The curve
defined  by \eqref{eq1} is an ellipse, when $\theta \in
[\piq,\pid[\,,$ a hyperbola, when $\theta \in ]\pid,\pitq[\,,$ the union
of two vertical lines, when $\theta = \pid\,,$ the diagonals
$\{x-y=0\}\cup \{x+y=0\}$, when $\theta = \pitq\,$. We only have to
analyze how the ellipses and the hyperbolas are situated within the
square.

\paragraph{Boundary points.}  At the boundary, for example on
$u= \pm 1$, \eqref{eq1} gives:
$$
  \cos \theta\,  ( 4v^2 - 1) + 3 \sin \theta  =0\,.
$$
 Depending on the value of $\theta$, we have no boundary point,
the zero set \eqref{eq1} is an ellipse contained in the open square;
one double boundary point, the zero set \eqref{eq1} is an ellipse,
contained in the closed square, which touches the boundary at two
points; or two boundary points, the zero set \eqref{eq1} is an
ellipse or a hyperbola meeting the boundary of the square at four
points.

\paragraph{Interior critical points.} We now look at the critical points of the
zero set of the function
$$
\Psi_{1,3} (u,v,\theta):=  \cos \theta\, ( 4v^2 - 1) + \sin \theta\,
(4u^2 - 1) \, .
$$
We get two equations:
$$
v \cos \theta =0\,,\, u \sin \theta =0\,.
$$
Except  for the two easy cases when $\cos \theta =0$ or $\sin \theta
=0$, which can be analyzed directly (product situation), we
immediately get that the only possible critical point is $(u,v)
=(0,0)$,  \ie \,  $(x,y) = ( \frac \pi 2,\frac \pi 2)$,  and that
this can only occur for $\cos \theta + \sin \theta =0\,$,  \ie \,
for $ \theta = \frac{3\pi}{4}\,$.

 The possible nodal patterns for an eigenfunction associated
with $\lambda_5$ all appear in Figure~\ref{F-pch}, second row.


This analysis shows rigorously that the number of nodal domains is
$2$, $3$ or $4$ as claimed in \cite{Pl}, and numerically observed in
Figure~\ref{F-pch}.  Observe that we have a rather complete
description of the situation by analyzing the points at the
boundary, and the critical points of the zero set inside the square.
When no critical point or no change of multiplicity is observed at
the boundary, the number of nodal domains remains constant. Hence,
the complete computation could be done by analyzing the ``critical''
values of $\theta$,  \ie  those for which there is a critical
point on the zero set of $\Psi_{1,3}$ in the interior, or a change
of multiplicity at the boundary, and one ``regular'' value of
$\theta$ in each non critical interval. To explore all possible
nodal patterns, and number of nodal domains, of the
eigenfunctions $\Phi^{\theta}_{1,3}, ~\theta \in [\piq,\pitq]\,$,
associated with the eigenvalue $\lambda_5$, it is consequently
sufficient to consider the values $\theta=\piq\,$,
$\theta=\arctan 3 \,$, $\theta =\frac \pi 2\,$, $\theta =
\frac{3\pi}{4}\,$.

\subsection{Second case: eigenvalue $\lambda_7$, or $(m,n)= (2,3)$.}

We look at the zeros of  $\Phi_{2,3}(x,y,\theta)$. We first observe
that
$$
\Phi_{2,3}(x,y,\theta)= \sin x \sin y \left( 2 \cos \theta \cos x
(\cos 2y + 2 \cos^2 y) + 2 \sin \theta \cos y  (\cos 2 x + 2 \cos^2
x)\right) \,.
$$

In  the coordinates \eqref{eq0}, this reads:
\begin{equation}
\Phi_{2,3}(x,y,\theta)=2 \sqrt{1-u^2} \sqrt{1-v^2} \left( u  \cos
\theta ( 4v^2 - 1) +  v \sin \theta (4u^2 - 1) \right)\,.
\end{equation}
We have to look at the solutions of
\begin{equation}\label{eqb1}
\Psi_{2,3}(u,v,\theta):= u ( 4v^2 - 1)\,  \cos \theta  +  v (4u^2 - 1)\,
\sin \theta  =0\,,
\end{equation}
inside $[-1,+1]\times [-1,+1]\,$.

\paragraph{Analysis at the boundary.}

The function $\Psi_{2,3}$ is anti-invariant under the change $(u,v)
\to (-u,-v)\,$. Changing $u$ into $-u$ amounts to changing $\theta$ in
$\pi-\theta$. Exchanging $(x,y)$ into $(y,x)$ amounts to changing
$\theta$ into $\frac{\pi}{2}-\theta\,$. This implies that it suffices
to consider the values $\theta \in [0,\frac{\pi}{2}]\,$, and the
boundaries $u=-1$ and $v=-1\,$. At the boundary $u=-1\,$, we get:
\begin{equation}\label{eqb1a}
-   \cos \theta \, ( 4v^2 - 1) +  3  v \sin \theta =0\,,
\end{equation}
with the condition that $v\in [-1,+1]\,$.

We note that the product of the roots is $-\frac 14\,$, and that there
are always two distinct solutions in $\mathbb R$. For $\theta =0\,$,
we have two solutions given by $v=\pm \frac 12\,$. When $\theta$
increases we still have two solutions in $[-1,1]$ till the largest
one is equal to $1$ (the other one being equal to $-\frac 14$). This
is obtained for $\theta = \frac \pi 4\,$. For $ \theta >
\frac{\pi}{4}\,$, there is only one negative solution in $[-1,1]\,$,
tending to zero as $\theta \ar \frac \pi 2\,$.

At the boundary $v=-1\,$, we have for $\theta =0\,$, $u=0$ as unique
solution. When $\theta$ increases, there is only one solution in
$[-1,1]$\,, till $\frac \pi 4$, where we get two solutions $u =- \frac
14$ and $u=1\,$. For this value of $\theta$, the zero set is given by
$(u+v) (4 uv -1)=0\,$. For $\theta \in ]\frac \pi 4,\frac \pi 2]\,$, we
have two solutions.

We conclude that the zero set of $\Psi_{2,3}$ always hits the
boundary at  six  points.

\paragraph{Critical points.}
We now look at the critical points of $\Psi_{2,3}$. We get two
equations:
\begin{equation} \label{eqb2}
( 4v^2 - 1)  \cos \theta  + 8u v \sin \theta =0\,,
\end{equation}
and
\begin{equation}\label{eqb3}
8  u  v \, \cos \theta  +   (4u^2 - 1) \, \sin \theta  =0\,.
\end{equation}
The critical points on the zero set of $\Psi_{2,3}$ are the common
solutions of \eqref{eqb1}, \eqref{eqb2}, and \eqref{eqb3}.

If $\cos \theta \sin \theta \neq 0$, we immediately obtain that
$u=v=0\,$, and these equations have no common solution. It follows
that the eigenfunctions associated with $\lambda_7$ have no interior
critical point on their nodal sets.

One can give the following expressions for the partial derivatives
of $\Psi_{2,3}$,
$$
\partial_u\Psi_{2,3}(u,v,\theta) = \frac{v}{u}(4u^2+1) \sin\theta
\text{~and~} \partial_v\Psi_{2,3}(u,v,\theta) = \frac{u}{v}(4v^2+1)
\cos\theta\,,
$$
for $u$, \resp $v$, different from $0$. Since a regular closed curve
contains points with vertical or horizontal tangents, it follows
that the zero set of $\Psi_{2,3}$ cannot contain any closed
component (necessarily without self-intersections, otherwise the
nodal set of $ \Phi_{2,3}$ would have a critical point). The
components of this zero set are lines joining two boundary points
which are decreasing from the left to the right. These lines cannot
intersect each other (for the same reason as before).

The possible nodal patterns for an eigenfunction associated
with $\lambda_7$ all appear in Figure~\ref{F-pch}, third row.

The number of nodal domains is  four (delimited by three non
intersecting lines) or  six   in the product case. Hence  the maximal
number of nodal domains is six.

\subsection{Third case : eigenvalue $\lambda_9$, or $(m,n)= (1,4)$}

We look at the zeros of $\Phi_{1,4} (\cdot,\cdot, \theta)$. Here  we
can write
\begin{equation*}
\Phi_{1,4}(x,y,\theta)= 4  \sin x \sin y\,  \Psi_{1,4} (u,v,\theta)
\end{equation*}
with
\begin{equation*}
\Psi_{1,4} (u,v,\theta) := \cos \theta\,  v (2v^2-1) + \sin \theta\,
u (2u^2-1) \,.
\end{equation*}
Hence, we have to analyze the equation
\begin{equation}\label{zeroc}
\cos \theta\,  v (2v^2-1) + \sin \theta\,  u (2u^2-1)=0\,.
\end{equation}

Notice that the functions $\Psi_{1,4}(u,v,\theta)$ are
anti-invariant under the symmetry $(u,v) \to (-u,-v)$, and that one
can reduce from $\theta \in [0,\pi[$ to the case $\theta \in
[0,\frac \pi 2]$ by making use of the symmetries with respect to the
lines $\{u=0\}$, $\{v=0\}$ and $\{u=v\}$.  One can even reduce
the analysis to $\theta \in [0,\piq]$ by changing $\theta$ into
$\pid-\theta$, and $(x,y)$ into $(y,x)\,$.

\paragraph{Boundary points.}  Due to the symmetries, the zero set
of $\Psi_{1,4}$ hits parallel boundaries at symmetrical points. For
$u=\pm 1$ these points are given by:
$$
v (2v^2-1) \pm  \tan \theta =0\,.
$$
If we start from $\theta=0\,$, we first have three zeroes,
corresponding to points at which the zero set of $\Psi_{1,4}$
arrives at the boundary: $0, \pm \frac{1}{\sqrt{2}}\,$.  Looking at
the derivative, we have a double point when $v=\pm \frac
{1}{\sqrt{6}}\,$, which corresponds to $\tan \theta =
\frac{\sqrt{2}}{3 \sqrt{3}}\,$. For larger values of $\theta$, we have
only one point till $\tan \theta =1\,$.

Hence, there are $3$, $2$, $1$ or $0$ solutions satisfying $v\in
[-1,+1]\,$. The analogous equation for $v = \pm 1$ appears with
$\cot\theta$ instead of $\tan\theta$, so that the boundary analysis
depends on the comparison of $|\tan \theta|$ with $\frac{\sqrt{2}}{3
\sqrt{3}}\,$, $1$ and $\frac{3 \sqrt{3}}{\sqrt{2}}\,$. When the points
disappear on $u=\pm 1$\,, they appear on $v=\pm 1\,$. Notice that the
value $\frac{\sqrt{2}}{3 \sqrt{3}}$ appears in Figure~\ref{F-pch} and
in Courant-Hilbert's book \cite{CH}, \S V.5.3, p. 302. Finally, the
maximal number of points along the boundary is  six, counting
multiplicities.

\begin{figure}[!bt]
\begin{center}
\includegraphics[width=14cm]{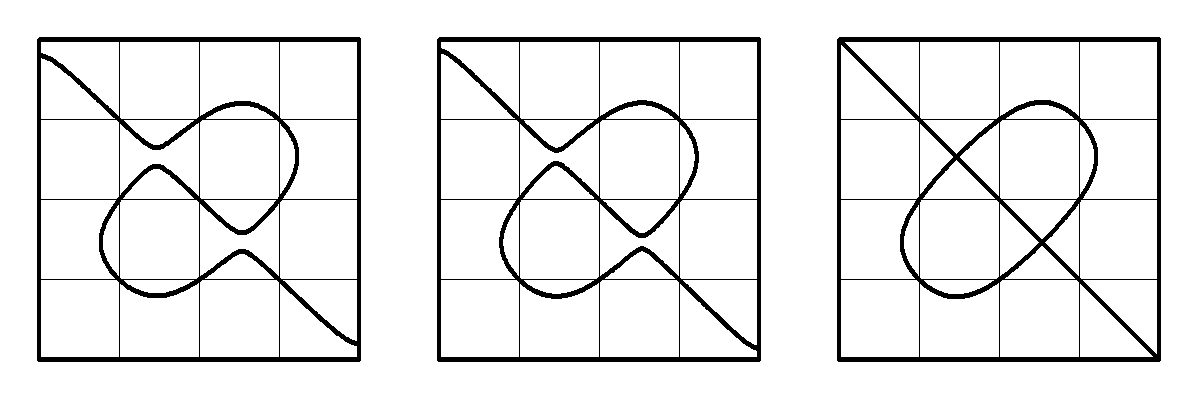}
\caption{Eigenvalue $\lambda_9$, deformation of the nodal set near
$\theta =\frac{\pi}{4}$}\label{fig-14-desing}
\end{center}
\end{figure}

\begin{figure}[!ht]
\begin{center}
\includegraphics[width=14cm] {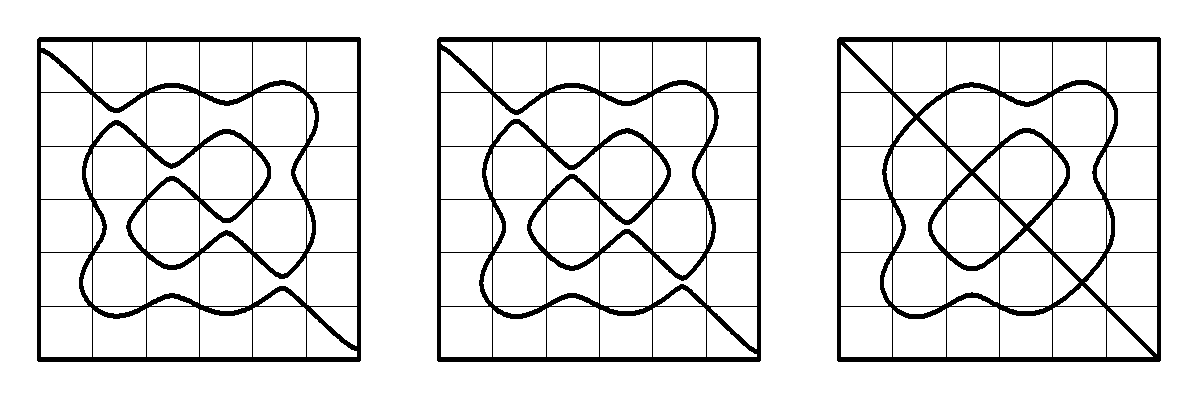}
\caption{Eigenvalue $\lambda_{23}$, deformation of the nodal set
near $\theta =\frac{\pi}{4}$}\label{fig-16-desing}
\end{center}
\end{figure}

\paragraph{Critical points.} The critical points of $\Psi_{1,4}$ satisfy:
\begin{equation}\label{eqc2}
\cos \theta \, (6 v^2-1)=0 \,,
\end{equation}
and
\begin{equation}\label{eqc3}
\sin \theta \, (6 u^2-1) =0\,.
\end{equation}
If we exclude the ``product'' case, the only critical points are
determined by $u^2=\frac 16$ and $v^2=\frac 16\,$. Plugging these
values in \eqref{zeroc}, we obtain that interior critical points on
the zero set of $\Psi_{1,4}$ can only appear when:
\begin{equation}\label{thetac}
\cos \theta\, \pm  \sin \theta\, =0\,.
\end{equation}
Hence, we only have to look at $\theta = \frac \pi 4$ and $\theta =
\frac{3\pi}{4}$. Because of symmetries, it suffices to consider the
case $\theta = \frac \pi 4$:
$$
\Psi_{1,4} (u,v,\frac \pi 4) := \frac {1}{\sqrt{2}}  ( v (2v^2-1) +
u (2u^2-1)) =   \frac {1}{\sqrt{2}}  ( u+v) ( 2 (u- \frac{v}{2}) ^2
+ \frac 32  v^2 -1)   \,.
$$
The zero set is the union of an ellipse contained in the square, and
 the anti-diagonal, with two intersection points. It follows
that the function $\Phi_{1,4}(x,y,\frac \pi 4)$ has four nodal
domains. Figure~\ref{fig-14-desing} shows the deformation of the
nodal set of $\Phi_{1,4}(x,y,\theta)$ for $\theta \le \frac{\pi}{4}$
close to $\frac{\pi}{4}\,$, as well as the grid $\{ \sin(4x)\sin(4y)=0 \}$. Figure~\ref{fig-16-desing} shows the deformation for $\Phi_{1,6}(x,y,\theta)$.

Let us summarize what we have so far obtained for the eigenfunctions
associated with $\lambda_9$.\vspace{-3mm}
\begin{itemize}
\item We have determined the  shape of the nodal set of
$\Phi_{1,4}$ when $\theta = \frac \pi 4$ or $\frac{3\pi}{4}\,$, and we
can easily see that these are the only cases  in  which the
interior part of the nodal set hits the boundary at the vertices.
\item When $\theta \not = \frac \pi 4$ or $\frac{3\pi}{4}\,$,
we have proved that the nodal set of $\Phi_{1,4}$ has no interior
critical point, hence no self-intersection, and that it hits the
boundary at $2$ or $6$ points, counting multiplicities.
\item We can also observe that all nodal sets must contain
the lattice points $(i\frac{\pi}{4},j\frac{\pi}{4})$ for $1 \le i,j
\le 3\,$, and that these points are always regular points of the
nodal sets. This implies, by energy considerations, that the nodal
sets cannot contain any closed component avoiding these lattice
points. The lattice points are indicated in Figure~\ref{F-pch}. They appear in Figures~\ref{fig-14-desing} and \ref{fig-14-desing} as the vertices of the grids.
\end{itemize}

 We still need to prove that all the possible nodal patterns
for an eigenfunction associated with the eigenvalue $\lambda_9$
appear in Figure~\ref{F-pch}, fourth row, and hence that the maximal
number of nodal domains is $4\,$, so that $\lambda_9$ is not Courant
sharp. When $\theta \in [0,\piq]\,$, this can be done by looking at
the intersections of the nodal sets with the horizontal lines
$\{y=\arccos(\pm \sqrt{\frac{1}{6}})\}$\,.  We leave the details to
the reader, and refer to Section~\ref{S-1R} for general arguments.

{ \textbf{Remark}}. Figure~\ref{F-pch} and \ref{fig-14-desing}
indicate that for some values of $\theta$, the function
$\Phi_{1,4}(x,y,\theta)$ has exactly two nodal domains. This
phenomenon  has been studied by Antonie Stern \cite{St} who
claims that for any $k \ge 2\,$, there exists an eigenfunction
associated with the Dirichlet eigenvalue $(1+4k^2)$ of the square
$[0,\pi]^2$, with exactly two nodal domains. In
Sections~\ref{HStern} -- \ref{S-1R}, we look at Stern's thesis more
carefully.


\section{The observations of A. Stern}\label{HStern}

The general topic of A.~Stern's thesis \cite{St} is the asymptotic
behaviour of eigenvalues and eigenfunctions. In Part I, she studies
the nodal sets of eigenfunctions of the  Laplacian in the
square with Dirichlet boundary condition, and on the $2$-sphere.
The eigenvalues are chosen to be positive, and listed in increasing
order, with multiplicities. In \cite{St1}, we propose extracts from
Stern's thesis, with annotations and highlighting to point out the
main results and ideas.

As we have seen in the previous sections, Pleijel's theorem
\cite{Pl} states that for a plane domain, there are only finitely
many Courant sharp Dirichlet eigenvalues.   For the Dirichlet
Laplacian in the square, A.~Stern claims \cite[tags E1, Q1]{St1}
that there are actually infinitely many eigenfunctions having
exactly two nodal domains,
\begin{quote}
[E1]\ldots Im eindimensionalen Fall wird nach den S\"atzen von
Sturm\footnote{Journal de Math\'ematiques, T.1, 1836, p. 106-186,
269-277, 375-444} das Intervall durch die Knoten der $\it n$ten
Eigenfunktion in $n$ Teilgebiete zerlegt. Dies Gesetz verliert seine
G\"ultigkeit bei mehrdimensionalen Eigenwertproblemen, \ldots es
l\"a{\ss}t sich beispielweise leicht zeigen, da{\ss} auf der Kugel
bei jedem Eigenwert die Gebietszahlen $2$ oder $3$ auftreten, und
da{\ss} bei Ordnung nach wachsenden Eigenwerten auch beim Quadrat
die Gebietszahl $2$ immer wieder vorkommt.\\

[Q1]\ldots Wir wollen nun zeigen, da{\ss} beim Quadrat die Gebietszahl
zwei immer wieder auftritt.
\end{quote}

 The second statement is mentioned in the book of
Courant-Hilbert \cite{CH}, \S VI.6, p. 455.

Pleijel's theorem has been generalized to surfaces by J. Peetre
\cite{Pe}, see also \cite{BeMe}.  As a consequence, only
finitely many eigenvalues of the sphere are Courant sharp. A.~Stern
claims  \cite[tags E1, K1, K2]{St1} that, for any $\ell \ge 2$,
there exists a spherical harmonic of degree $\ell$ with exactly
three nodal domains when $\ell$ is odd, \resp with exactly two nodal
domains when $\ell$ is even,
\begin{quote}
[K1] \ldots Zun\"achst wollen wir zeigen, da{\ss} es zu jedem Eigenwert
Eigenfunktionen gibt, deren Nullinien die Kugelfl\"ache nur in zwei
oder drei Gebiete teilen.\\

[K2] \ldots  Die Gebietszahl zwei tritt somit bei allen
Eigenwerten $\lambda_n = (2r+1)(2r+2)\qquad r=,1,2,\cdots $ auf;
ebenso wollen wir jetzt zeigen, da{\ss} die Gebietszahl drei bei allen
Eigenwerten
$$
\lambda_n =2r  (2r+1)\qquad r=1,2,\cdots
$$
immer wieder vorkommt.
\end{quote}
These two statements are usually attributed to H.~Lewy
\cite{Lew}.\medskip

In this paper, we shall only deal with the case of the square,
leaving the case of the sphere for \cite{BeHe2}. First, we
quote the main statements made by A.~Stern \cite[tags Q1-Q3]{St1},
and summarize them in Theorem~\ref{stern-T}.
\begin{quote}
[Q2]\ldots Wir betrachten die Eigenwerte
$$ \lambda_n= \lambda_{2r,1} = 4r^2 +1\,,\, r=1,2,\dots$$
und die Knotenlinie der zugeh\"orige Eigenfunktion
$$
u_{2r,1} + u_{1,2r} =0\,,
$$
f\"ur die sich, wie leicht mittels graphischer Bilder nachgewiesen
werden kann, die Figur~7 ergibt.\\

[Q3]\ldots La{\ss}en wir nur $\mu$ von $\mu=1$ aus abnehmen, so l\"osen
sich die Doppelpunkte der Knotenlinie alle gleichzeitig und im
gleichem Sinne auf, und es ergibt sich die Figur~8. Da die
Knotenlinie aus einem Doppelpunktlosen Zuge besteht, teilt sich das
Quadrat in zwei Gebiete und zwar geschieht  dies f\"ur alle Werte
$r=1,2,\dots\,,$ also Eigenwerte $\lambda_{n} =
\lambda_{2r,1}=4r^2+1\,$.
\end{quote}

\begin{figure}[!ht]
\begin{center}
\includegraphics[width=12cm]{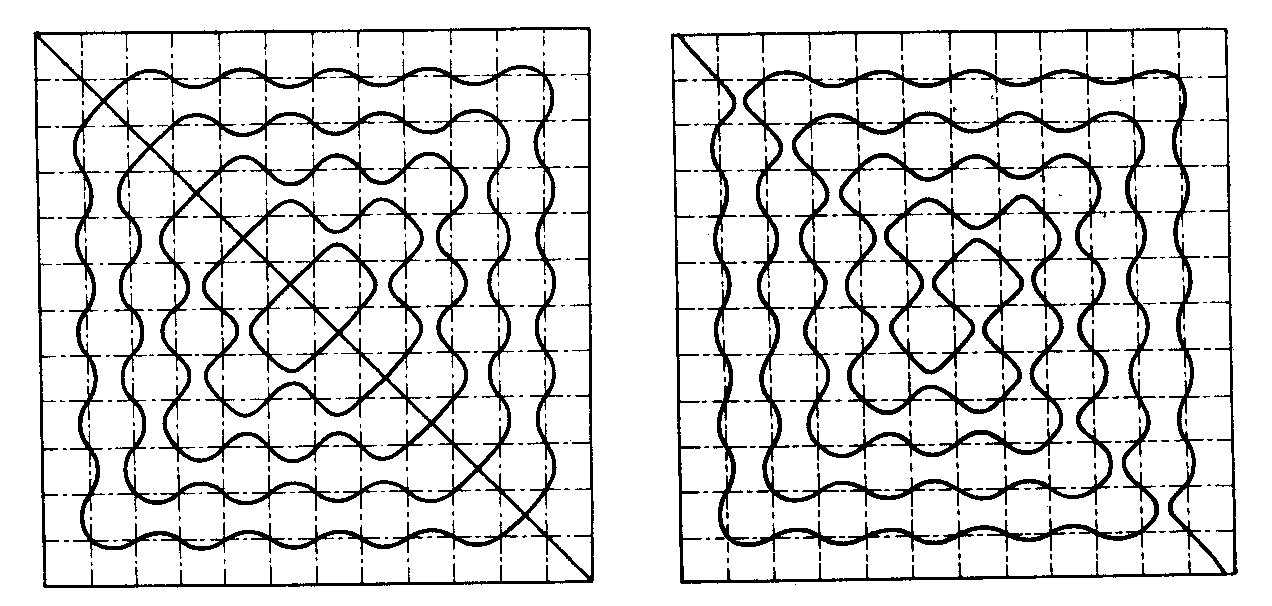}
\caption{Case $r=6$, nodal sets for $\theta=\frac \pi 4$ and
$\theta$ close to $\frac \pi 4$ (reproduced from \cite{CHd})}\label{stern-0}
\end{center}
\end{figure}

\begin{theorem}\label{stern-T} For any $r\in \mathbb N$, consider the family
$\Phi_{1,2r}(x,y,\theta)$ of eigenfunctions of the Laplacian in the
square $[0,\pi]^2$, associated with the Dirichlet eigenvalue $1+
4r^2$,
$$\Phi^{\theta}_{1,2r}(x,y) := \Phi_{1,2r}(x,y,\theta) :=\cos \theta \sin x \sin
(2r  y) + \sin \theta \sin (2r x) \sin y\,.$$ Then,\vspace{-3mm}
\begin{enumerate}
\item  for $\theta =\frac \pi 4$, the nodal pattern of $\Phi$ is as
shown in Figure~\ref{stern-0}, left, \cite[Fig.~7]{St};
\item for  $\theta <\frac \pi 4\,$, and $\theta$ sufficiently close to
$\frac \pi 4$, the double points all disappear at the same time and
in a similar manner (`im gleichem Sinne') as in
Figure~\ref{stern-0}, right, \cite[Fig.~8]{St}. The nodal set consists of a connected line (`aus
einem Zuge') with no double point. It divides the square into two
domains.
\end{enumerate}
\end{theorem}

{\textbf{Remark}}. Although this is not stated explicitly,  one
can infer from Stern's thesis, (i) that the eigenfunction $
\Phi_{1,2r}(x,y,\frac \pi 4)$ has $2r$ nodal domains and $(2r-2)$
double points, and (ii) that for $\theta$ close to and different
from $\frac \pi 4$, the nodal set of $ \Phi_{1,2r}(x,y,\theta)$
consists of the boundary of the square and a \emph{connected} simple
curve from one point of the boundary to a symmetric point.  This
curve divides the domain into two connected components.

\medskip

\begin{figure}[!ht]
\begin{center}
\includegraphics[width=4cm]{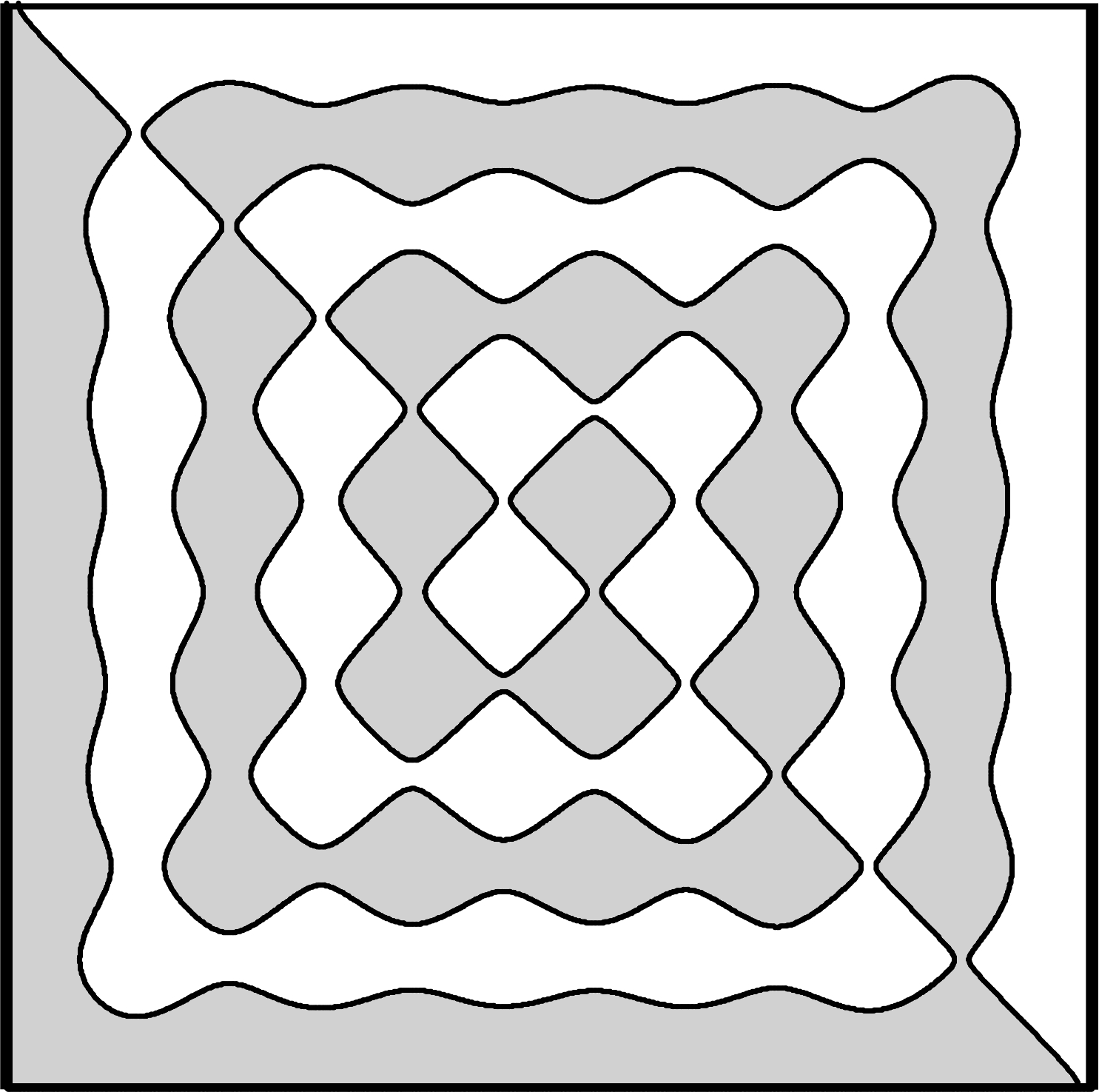}
\caption{Nodal domains, courtesy Virginie Bonnaillie-No\"{e}l
\cite{vbn}}\label{vbn}
\end{center}
\end{figure}

A.~Stern states two simple properties which play a key role in
the proofs \cite[tags I1, I2]{St1}. These statements are formalized
in Properties~\ref{stern-P1} below.
\begin{quote}
[I1]\ldots Um den typischen Verlauf der Knotenlinie zu bestimmen,
haben wir \"ahnliche Anhaltspunkte wie auf der Kugelfl\"ache. Legen
wir die Knotenliniensysteme von $u_{\ell,m}$ ($\ell -1$ Parallelen
zur $y$-Achse, $m-1$ zur $x$-Achse) und $u_{m,\ell}$ ($m -1$
Parallelen zur $y$-Achse, $\ell-1$ zur $x$-Achse) \"ubereinander,
 so kann f\"ur $\mu >0$ ($<0$) die Knotenlinie nur in den
Gebieten verlaufen, in denen beide Funktionen verschiedenes
(gleiches) Vorzeichen haben. \\

[I2]\ldots Weiter m\"ussen alle zum Eigen\-wert $\lambda_{n,m}$
geh\"origen Knotenlinien durch Schnitt\-punkte der Liniensysteme
$u_{\ell,m}=0$ and $u_{m,\ell}=0$, also durch $(\ell-1)^2 + (m-1)^2$
feste Punkte hindurchgehen \ldots
\end{quote}

\begin{properties}\label{stern-P1}
Let $\phi$ and $\psi$ be two linearly independent eigenfunctions
associated with the same eigenvalue for the square $\Sg$. Let $\mu$
be a real parameter, and consider the family of eigenfunctions
$\phi_{\mu} = \psi + \mu \phi$. Let $N(\phi)$ denote the nodal set
of the eigenfunction $\phi$.
\begin{enumerate}
    \item Consider the domains in $\Sg \setminus N(\phi) \cup
    N(\psi)$ in which $\mu \,\phi \,\psi > 0$ and hatch them
    \footnote{`schraffieren', see \cite[tag I1]{St1}, in the
    spherical case.}.
    Then the nodal set $N(\phi_{\mu})$ avoids the hatched domains.
    \item The points in $N(\phi) \cap N(\psi)$ belong to the nodal
    set $N(\phi_{\mu})$ for all $\mu\,$.
\end{enumerate}
\end{properties}

\begin{property}\label{stern-P2}
The nodal set $N(\phi_{\mu})$ depends continuously on $\mu\,$.
\end{property}

\textbf{Remark}. As a matter of fact, A. Stern uses
Property~\ref{stern-P1} in both cases (square and sphere), and only
mentions Property~\ref{stern-P2} in the case of the sphere. She says
nothing on the proof of this  second  property which is more or
less clear near regular points, but not so clear near multiple
points. H.~Lewy provides a full proof in the case of the sphere
\cite[Lemmas~2-4]{Lew}.

Finally, A.~Stern mentions  that  she uses a graphical method
(`mittels graphischer Bilder' and `unter Zuhilfenahme graphischer
Bilder' \cite[tags Q2, Q4]{St1}) which may have been classical at
her time, and could explain the amazing quality of her pictures. On
this occasion, she mentions a useful idea in her \S~I.3, namely
looking at the intersections of the nodal set $N(\phi_{\mu})$ with
horizontal or vertical lines.\medskip

All in all, the arguments given by A.~Stern seem rather sketchy to
us, and we do not think that they are quite sufficient to conclude
the proof of Theorem~\ref{stern-T}.

In our opinion, taking care of the following items is missing in
Stern's thesis.\vspace{-3mm}
\begin{enumerate}
\item Complete determination of the multiple points of
$N(\Phi^{\frac{\pi}{4}})$ ;
\item Absence of multiple points in $N(\Phi^{\theta})$,
when $\theta$ is different from $\frac{\pi}{4}$, and close to
$\frac{\pi}{4}$ ;
\item Connectedness of the nodal set $N(\Phi^{\theta})$, or
why there are no other components, \eg closed inner components, in
the nodal set.
\end{enumerate}

The aim of this paper is to complete the proofs of A.~Stern in the
case of the square.

\textbf{Remark}. In \cite{Lew}, H.~Lewy  gives a complete
proof of Stern's results \footnote{ A.~Stern and H.~Lewy were both
students of R.~Courant at about the same time, 1925. H.~Lewy does
however not refer to A.~Stern's Thesis in his paper. We refer to
\cite{BeHe2} for a further discussion.} in the case of the sphere.
 In the unpublished preprint \cite{GSP}, the authors provide
partial answers to the above items in the case of the square.

 The key steps to better understand the possible nodal patterns
for the eigenvalues $ 1 + 4r^2$ (and other eigenvalues as well),
and to answer the above items, are the following. \vspace{-3mm}
\begin{itemize}
    \item Subsection~\ref{SS-1RCZ}, in which we study the points which are both
    zeroes and critical points of the functions $\Phi^{\theta}_{1,R}$.
    \item Subsection~\ref{SS-1RQ}, in which we study the possible local nodal
    patterns of the functions  $\Phi^{\theta}_{1,R}$.
    \item Subsection~\ref{SS-1RZ}, in which we determine the nodal sets
    of the functions  $\Phi^{\theta}_{1,R}$ for $\theta =
    \frac{\pi}{4}$ or $\frac{3\pi}{4}$.
\end{itemize}\vspace{-3mm}
In the subsequent subsections, we study the deformation of the nodal
set of $\Phi^{\theta}_{1,R}$ when $\theta$ varies close to
$\frac{\pi}{4}$ or $\frac{3\pi}{4}$, and conclude the proof of
Theorem~\ref{stern-T}.  As a matter of fact, our approach gives
the maximal interval in which the nodal set of $\Phi^{\theta}_{1,R}$
remains connected, without critical points, see
Lemma~\ref{1RD-L1}~(i).\medskip

\textbf{Sketch of the proof of Theorem~\ref{stern-T}}. Consider the
eigenvalue $\hat{\lambda}_{1,R}:= 1+R^2$ for the square $\Sg :=
]0,\pi[^2$ with Dirichlet boundary conditions, and consider the
eigenfunction
$$
\Phi^{\theta}(x,y) := \Phi(x,y,\theta) := \cos\theta \, \sin x
\,\sin(Ry) + \sin\theta \, \sin(Rx)\,\sin y \,, ~~\theta \in
[0,\pi[\,.
$$

Let us introduce the $Q$-squares,
$$Q_{i,j} := ]\frac{i\pi}{R},\frac{(i+1)\pi}{R}[ \times
]\frac{j\pi}{R},\frac{(j+1)\pi}{R}[\,, \text{~for~} 0 \le i,j \le
R-1\,,
$$
and the lattice,
$$
\Lg := \left\{(\frac{i\pi}{R},\frac{j\pi}{R}) ~|~ 1 \le i,j\le
R-1\right\}\,.
$$

The basic idea is to start from the analysis of a given nodal set,
\eg from the nodal set $N(\Phi^{\frac{\pi}{4}})$, and  then to
use some kind of perturbation argument.\medskip

Here are the key points.\vspace{-3mm}
\begin{enumerate}
\item Use Property~\ref{stern-P1}: Assertion (i) defines checkerboards by
$Q$-squares (depending on the sign of $\cos\theta$),  whose grey
squares do not contain any nodal point of $\Phi^{\theta}$. Assertion
(ii) says that the lattice $\Lg$ is contained in the nodal set
$N(\Phi^{\theta})$ for all $\theta$.
\item Determine the possible \emph{critical zeroes} of the eigenfunction
$\Phi^{\theta}$, \ie the zeroes which are also critical points,
both in the interior of the square or on the boundary. They indeed
correspond to multiple points in the nodal set. Note that the points
in $\Lg$ are not critical zeroes, see Subsection~\ref{SS-1RCZ}.
\item Determine whether critical zeroes are degenerate or not and
their order when they are degenerate.
\item Determine how critical zeroes appear or disappear when
$\theta$ varies, and how the nodal set $N(\Phi^{\theta})$ evolves.
For this purpose, make a local analysis in the square $Q_{i,j}\,$,
depending on whether it is contained in $\Sg$ or touches the
boundary, see Subsection~\ref{SS-1RQ}.
\item Determine the nodal sets of the eigenfunctions $Z_{\pm}$ associated
with the eigenvalue $\hat{\lambda}_{1,R}$. For this purpose,
determine precisely the critical zeroes of $\Phi^{\theta}$ for
$\theta = \frac{\pi}{4}$ and $\frac{3\pi}{4}$, and prove a
separation lemma in the $Q_{i,j}$ to determine whether the medians
of this $Q$-square meet the nodal set of $\Phi^{\theta}$ when
$\theta = \frac{\pi}{4}$ or $\frac{3\pi}{4}$, see
Subsection~\ref{SS-1RQ}.
\item Prove that the nodal set $N(\Phi^{\theta})$ does not contain
any closed component.
\end{enumerate}

Take $R=2r$ and $0 < \frac{\pi}{4} - \theta \ll 1\,$. Using the above
analysis one can actually give a complete proof of
Theorem~\ref{stern-T}. The analysis of the local possible nodal
patterns shows that the nodal set $N(\Phi^{\frac{\pi}{4}})$ is
indeed as stated. For $0 < \frac{\pi}{4} - \theta \ll 1\,$, the
eigenfunction $\Phi^{\theta}$ has no critical zero in $\Sg$, and
exactly two critical zeroes on the boundary, symmetric with respect
to the center of the square. This proves in particular that the
critical zeroes of $\Phi^{\frac{\pi}{4}}$ all disappear at once when
$\theta$ changes, $\theta < \frac{\pi}{4}\,$. Starting from one of the
critical zeroes on the boundary, and using the above analysis, one
can actually follow a connected nodal simple curve passing through
all the points in $\Lg$ and going from on of the critical zeroes on
the boundary to the second one. To finish the proof it suffices to
show that there are no other component of $N(\Phi^{\theta})$ in
$\Sg$.


\section{Notation and definitions. General properties of the nodal sets.}\label{S-G}

\subsection{Notation and definitions, I.}\label{SS-GNot1}

Let $\Sg$ be the open square $]0,\pi[^2$ in the plane. We denote by
$\partial \Sg$ boundary, by $\Dg_{+}$ the diagonal, by $\Dg_{-}$ the
anti-diagonal, and by $O := (\frac{\pi}{2},\frac{\pi}{2})$ the
center of $\Sg$.\medskip

Let $\Phi$ be an eigenfunction for the Dirichlet Laplacian in $\Sg$.
We let
\begin{equation}\label{G-not1}
N(\Phi):= \big\{(x,y) \in \overline{\Sg} ~|~ \Phi(x,y) = 0\big\}
\end{equation}
denote the nodal set of $\Phi$, and
\begin{equation}\label{G-not1i}
N_i(\Phi) := N(\Phi) \cap \Sg
\end{equation}
denote the interior part of $N(\Phi)$. \medskip

Given two integers $m, n \ge 1$, we consider the one-parameter
family of eigenfunctions,
\begin{equation}\label{G-1}
\Phi^{\theta}_{m,n}(x,y) := \Phi_{m,n}(x,y,\theta) := \cos\theta
\sin(mx) \sin(ny) + \sin\theta \sin(nx) \sin(my)\,,
\end{equation}
with $x,y \in [0,\pi]$ and $\theta \in [0,\pi[\,$. \\
Unless necessary,
we skip the index $(m,n)$. These eigenfunctions are associated with
the eigenvalue
\begin{equation}\label{G-1vp}
\hat{\lambda}_{m,n} := m^2+n^2\,.
\end{equation}

The following eigenfunctions are of parti\-cu\-lar interest.
\begin{equation}\label{G-2}
\begin{array}{lll}
X := \Phi^{0}\,, & Y  := \Phi^{\frac{\pi}{2}}\,,\\
Z_{+} := \Phi^{\frac{\pi}{4}}\,, &Z_{-}  := \Phi^{\frac{3\pi}{4}}\,.\\
\end{array}
\end{equation}

Denote by
\begin{equation}\label{G-2f}
\Lg := N_i(X) \cap N_i(Y)\,,
\end{equation}
the set of zeroes which are  common to all eigenfunctions $\Phi^{\theta},
\theta \in [0,\pi[\,$.

\begin{definition}\label{S-D1}
A \emph{critical zero} of $\Phi$ is a point $(x,y) \in
\overline{\Sg}$ such that both $\Phi$ and $\nabla \Phi$ vanish at $(x,y)$.
\end{definition}

\subsection{General properties of nodal sets.}\label{SS-GProp}

Although  stated  in the case of the square, the following
properties are quite general (see \cite{Be} and references therein)  for eigenfunctions of the Dirichlet realization of the Laplacian in a regular domain of $\mathbb R^2$.

\begin{properties}\label{G-P1}
Let $(x,y)$ be a point in $\Sg$ (an interior point).\vspace{-3mm}
\begin{enumerate}
    \item A non-zero eigenfunction $\Phi$ cannot vanish to infinite order at $(x,y)$.
    \item If the non-zero eigenfunction $\Phi$ vanishes at $(x,y)$, then the leading
    part of its Taylor expansion at $(x,y)$ is a harmonic homogeneous polynomial.
    \item If the point $(x,y)$ is a critical zero of the eigenfunction $\Phi$, then
    the nodal set $N(\Phi)$ at the point $(x,y)$ consists of finitely many
    regular arcs which form an equi-angular system.
    \item The nodal set can only have self-intersections at critical
    zeroes, and the number of arcs which meet at a
    self-intersection is determined by the order of vanishing of the
    eigenfunction. Nodal curves cannot meet tangentially.
    \item The nodal set cannot have an end point in the interior of
    $\Sg$, and consists of finitely many analytic arcs.
    \item Let the eigenfunction $\Phi$ be associated with the
    eigenvalue $\lambda$. If $\omega$ is a \emph{nodal domain},  \ie\, a
    connected component of $\Sg \setminus N(\Phi)$, then the first
    Dirichlet eigenvalue of $\omega$ is equal to $\lambda$.
    \item Similar properties hold at boundary points, in particular
    property (iii).

\end{enumerate}
\end{properties}

\textbf{Remark}. Since the eigenfunctions of $\Sg$ are
defined over the whole plane, the analysis of the critical zeroes at
interior points easily extends to the boundary.\medskip

\begin{properties}\label{G-P2}
Let $\Phi$ be an eigenfunction $\Phi_{m,n}^{\theta}$ of the square
$\Sg = ]0,\pi[^2\,$, with $\theta \in [0,\pi[\,$.\vspace{-3mm}
\begin{enumerate}
    \item For $\theta \not = \frac{\pi}{2}\,$, the nodal set $N(\Phi)$
    satisfies
    \begin{equation}\label{G-2a}
    \Lg \cup \partial\Sg \subset N(\Phi) \subset \Lg \cup
    \partial\Sg \cup \big\{ (x,y)\in [0,\pi]^2 ~|~ \cos\theta\,
    X(x,y)\, Y(x,y) < 0\big\}.
    \end{equation}
    \item If $\mathrm{gcd}(m,n) = 1$, then all the points in $\Lg$ are regular
    points of the nodal set.
    \item The nodal set $N(\Phi)$ can only hit the boundary of the
    square at critical zeroes (either in the interior of the edges
    or at the vertices).
    \item The nodal set $N(\Phi)$ can only pass from
    one connected component of the set
    $$\Wg_{m,n}^{\theta} := \big\{ (x,y)\in [0,\pi]^2 ~|~ \cos\theta\,
    X(x,y)\, Y(x,y) < 0\big\}$$
    to another through one of the points in $\Lg$.
    \item No closed connected component of $N(\Phi)$ can be contained in the closure
    of one of the connected components of $\Wg_{m,n}^{\theta}\,$. Equivalently, any
    connected component of $N_i(\Phi)$ must contain at least one point in $\Lg$.
\end{enumerate}
\end{properties}

 \textbf{Proof}. (i)  We have $\sin\theta > 0\,$, so that for
$\cos\theta\, X(x,y)\, Y(x,y) \ge 0$ the function $\Phi$ is either
positive or negative, it cannot vanish unless $(x,y) \in \Lg$. (ii)
Follows by direct analysis. (iii) Follows from Property~\ref{G-P1}.
(iv) Clear. (v) Any connected component of $N(\Phi)$ which does not
meet $\Lg$ would be strictly strictly contained in one of the nodal
domains of the eigenfunctions $X$ or $Y$, a contradiction with
Property~\ref{G-P1}~(vi). \hfill $\square$

\begin{figure}[!ht]
\begin{center}
\includegraphics[width=14cm]{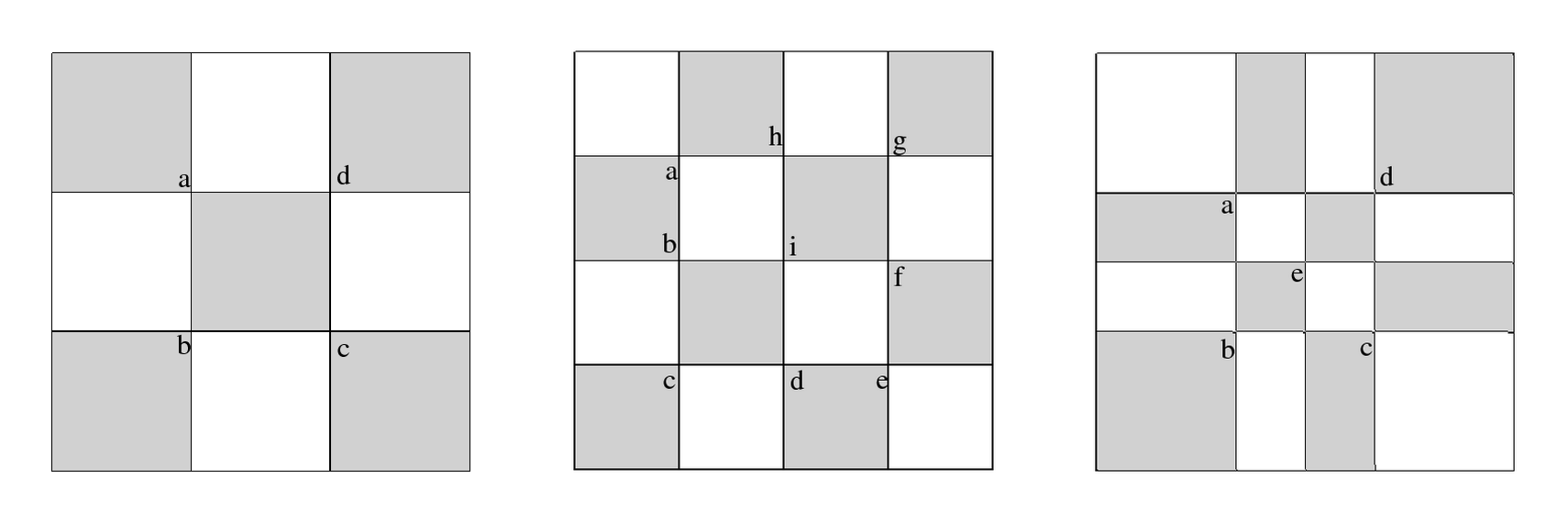}
\caption{Checkerboards for eigenvalues $\hat{\lambda}_{1,3}$\,,
$\hat{\lambda}_{1,4}$ and $\hat{\lambda}_{2,3}$}\label{FP-1}
\end{center}
\end{figure}

Figure~\ref{FP-1} illustrates property (i) when $(m,n) = (1,3),
(1,4)$ or $(2,3)$. When $\cos\theta > 0$, the nodal set is contained
in the white sub-squares; when $\cos\theta < 0$ it is contained in
the grey sub-squares. The points in $\Lg$ are the points labelled
$a, b, \ldots$ in the figures.

\subsection{Notation and definitions, {II}.}\label{SS-GNot2}

We now consider the case of the eigenvalue $\hat{\lambda}_{1,R} =
1+R^2$, for some integer $R\ge 1\,$.\\
 We introduce\\[5pt]
\noib the numbers
\begin{equation}\label{1RNot-2}
p_i := i \frac{\pi}{R}\,, \text{~for~} 0 \le i \le R\,,
\end{equation}
\begin{equation}\label{1RNot-3}
m_i := (i+\frac{1}{2}) \frac{\pi}{R}\,, \text{~for~} 0 \le i \le R-1\,,
\end{equation}
\noib the collection of squares
\begin{equation}\label{1RNot-4}
Q_{i,j} := ~ ]p_i,p_{i+1}[ \times ]p_j,p_{j+1}[\,, \text{~for~~} 0 \le
i,j \le R-1\,,
\end{equation}
whose centers are the points $(m_i,m_j)$,\\
\noib the lattice
\begin{equation}\label{1RNot-5}
\Lg := \left\{ (p_i,p_j) ~|~ 1 \le i,j \le R-1\right\}.
\end{equation}\bigskip

\textbf{Coloring the squares}. Assume that $\theta \not = 0$ and
$\frac{\pi}{2}$. If $(-1)^{i+j} \cos\theta < 0$, we color the square
$Q_{i,j}$ in white, otherwise we color it in grey. The collection of
squares $\{Q_{i,j}\}$ becomes a grey/white checkerboard (which
depends on $R$ and on the sign of $\cos\theta$). Depending on the
sign of $\cos\theta$, the white part of the checkerboard is given
by,
\begin{equation}\label{1RNot-6}
\begin{split}
\Wg(+) &:= \bigcup_{(-1)^{i+j}=-1} Q_{i,j}\,, \text{~when~} \cos\theta > 0\,,\\
\Wg(-) &:= \bigcup_{(-1)^{i+j}=1} Q_{i,j}\,,  \text{~when~} \cos\theta < 0\,.\\
\end{split}
\end{equation}

For the eigenfunction $\Phi^{\theta}$, we have,
\begin{equation}\label{1RNot-7}
\Lg \cup \partial \Sg \subset N(\Phi^{\theta}) \subset \Wg(\pm) \cup
\Lg \cup
\partial \Sg\,,
\end{equation}
if $(\pm\,\cos\theta > 0)$.\medskip

\textbf{Remark}. Observe that the squares $Q_{i,j}$ are open,
the sets $\Wg(\pm)$ do not contain the segments $\{x=p_i\}\cap \Sg$
 and $\{y=p_j\}\cap \Sg\,$.\medskip

\begin{figure}[!ht]
\begin{center}
\includegraphics[width=12cm]{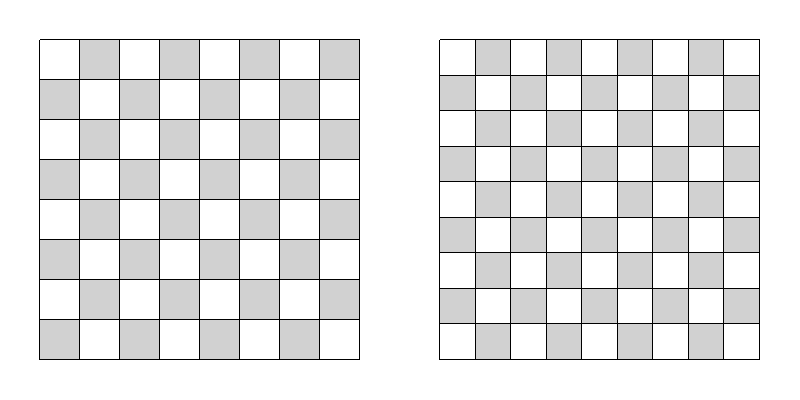}\vspace{-3mm}
\caption{Checkerboards $\Wg(+)$ for the eigenvalue
$\hat{\lambda}_{1,8}$, and $\Wg(-)$ for the eigenvalue
$\hat{\lambda}_{1,9}\,$} \label{1RNot-F1}
\end{center}
\end{figure}

Figure~\ref{1RNot-F1} shows the checkerboards for the eigenvalue
$\hat{\lambda}_{1,8}$, when $\cos\theta > 0\,$, \resp for
$\hat{\lambda}_{1,9}$, when $\cos\theta < 0\,$.

To describe the global aspect of the nodal sets, we will also use
the following squares.

Denote by
\begin{equation}\label{1RNot-8}
r := [\frac{R}{2}]\,,
\end{equation}
the integer part of $R/2\,$. For $0 \le i \le r\,$, define the square
\begin{equation}\label{1RNot-9}
\Sg_i := ~ ]p_i,p_{R-i}[ \times ]p_i,p_{R-i}[\,.
\end{equation}
With this notation, we have
$$\Sg_r \subset \Sg_{r-1}\subset \cdots \subset \Sg_0 = \Sg.$$
Furthermore, when $R=2r$, $\Sg_{r-1} = ]p_{r-1},p_{r+1}[^2$ consists
of four $Q$-squares, while $\Sg_r$ is empty; when $R=2r+1$,
$\Sg_{r}$ is a single $Q$-square. All these squares have the same
center $O = (\pi/2,\pi/2)$.


\section{Eigenfunctions associated with the eigenvalue $\hat{\lambda}_{1,R}$}\label{S-1R}

In this section, we consider the eigenfunctions associated with the
eigenvalue $\hat{\lambda}_{1,R}$, for an integer $R \ge 1$. More
precisely, we consider the $1$-parameter family of eigenfunctions,
\begin{equation}\label{1R-1}
\Phi^{\theta}(x,y) := \Phi(x,y,\theta) := \cos\theta \, \sin x \,
\sin(Ry) + \sin\theta \sin(Rx) \, \sin y \,,
\end{equation}
where $x, y \in [0,\pi]^2$ and $\theta \in [0,\pi[\,$.

This eigenfunction can be written as
\begin{equation}\label{1RChe-2}
\Phi(x,y,\theta) = \sin x \, \sin y \, \phi(x,y,\theta)\,,
\end{equation}
 with
\begin{equation}\label{1RChe-3}
\phi(x,y,\theta) := \cos\theta \, U_{R-1}(\cos y) + \sin\theta \,
U_{R-1}(\cos x)\,,
\end{equation}
where $U_n(t)$ is the $n$-th Chebyshev polynomial of second type
defined by the relation,
\begin{equation}\label{1RChe-1}
\sin t \; U_n(\cos t) := \sin\big( (n+1)t \big)\,.
\end{equation}

\subsection{Chebyshev polynomials and special values of $\theta$}\label{SS-1RChe}

In this section, we list some properties of the Chebyshev
polynomials to be used later on.

\begin{properties}\label{1RChe-P1}
For $ R \in \nzb\,$, the Chebyshev polynomial $U_{R-1}(t)$ has the
following properties.\vspace{-3mm}
\begin{enumerate}
\item The polynomial $U_{R-1}$ has degree $R-1$ and the same parity as
$R-1$. Its zeroes are the points $\cos p_j, 1 \le j \le R-1$, see
\eqref{1RNot-2}. Furthermore, $U_{R-1}(1)=R$, $U_{R-1}(-1) =
(-1)^{R-1}R\,$, and $-R \le U_{R-1}(t) \le R$ for all $t\in [-1,1]\,$.
\item The polynomial $U_{R-1}'$ has exactly $R-2$ simple zeroes
$\cos q_j, 1\le j \le R-2$, with $q_j \in ]p_j,p_{j+1}[$.
\item When $R$ is even, $R=2r$, the values $q_j$ satisfy,
\begin{equation}\label{1RChe-4e}
\begin{split}
0 < q_1 < q_2 \cdots < q_{r-1} & < \frac{\pi}{2} < q_r < \cdots <
q_{2r-2} < \pi\,,\\
q_{2r-1-j} & = \pi - q_{j}\,, ~1\le j \le r-1\,.\\
\end{split}
\end{equation}
\item When $R$ is odd, $R=2r+1\,$, the values $q_j$ satisfy,
\begin{equation}\label{1RChe-4o}
\begin{split}
0 < q_1 < q_2 \cdots < q_{r-1} < q_r &= \frac{\pi}{2} < q_{r+1} <
\cdots < q_{2r-1} < \pi\,,\\
q_{2r-j} & = \pi - q_{j}, ~1\le j \le r-1\,.\\
\end{split}
\end{equation}
\item Let $M_j := U_{R-1}(\cos q_j)$, for $1\le j \le R-2\,$, denote the local extrema of
$U_{R-1}$. Then,
\begin{equation}\label{1RChe-5}
\begin{split}
& (-1)^{j} M_j > 0 \text{~and~} (-1)^{j} U_{R-1}(\cos t) > 0 \text{~in~} ]p_j,p_{j+1}[\,,\\
& (-1)^{j+1} \big( U_{R-1}(\cos t) - M_j \big) \ge 0 \text{~in~}
]p_j,p_{j+1}[\,.\\
\end{split}
\end{equation}
\end{enumerate}
\end{properties}

\begin{figure}[!ht]
\begin{center}
\includegraphics[width=12cm]{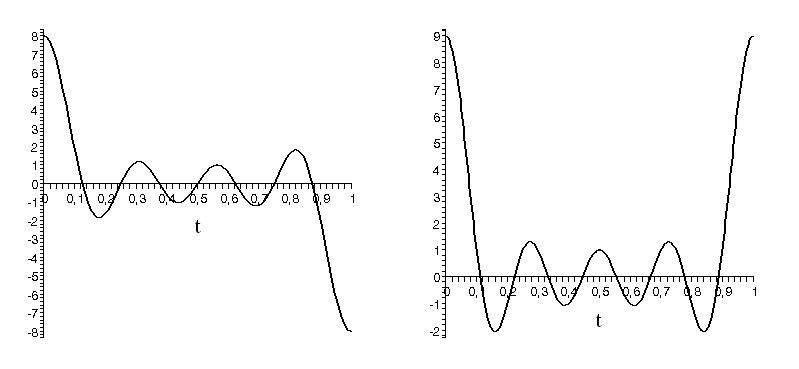}
\caption{Functions $U_7(\cos(\pi t))$ and $U_8(\cos (\pi t))$}
\label{1RChe-F1}
\end{center}
\end{figure}

 \textbf{Proof}. The above properties are easy to prove, and
illustrated by the graph of the function $t \to U_{R-1}(\cos t)$ in
the interval $[0,\pi]$, see Figure~\ref{1RChe-F1}, for the cases
$R=8$ and $R=9$. \hfill $\square$ \bigskip

\textbf{Special values of the parameter $\theta$}. We shall now
associate some special values of the parameter $\theta$ with the
zeroes
\begin{equation}\label{1-RChe-6}
\Qg := \left\{ q_j ~|~ 1\le j \le R-2 \right\}
\end{equation}
of the function $t \to U_{R-1}'(\cos t)$\,. As we shall see later on,
they are related to changes in the nodal patterns of the
eigenfunctions $\Phi^{\theta}$ when $\theta$ varies from $0$ to
$\pi$.\medskip

The values of $\theta$ to be introduced below are well defined
because the polynomial $U_{R-1}$ does not vanish at the points $\cos
q_k$, $1 \le k \le R-2\,$. These values of $\theta$ will clearly
depend on $R$, although we do not indicate the dependence in the
notations.
\medskip

\noib For $1 \le i,j \le R-2$, define $\theta(q_i,q_j)$, alias
$\theta_{i,j}$, to be the unique angle in the interval  $[0,\pi[$
such that
\begin{equation}\label{1-RChe-6a}
\cos\theta_{i,j} \, U_{R-1}(\cos q_j) + \sin\theta_{i,j} \,
U_{R-1}(\cos q_i) = 0\,.
\end{equation}
Let $\Tg_o$ denote the corresponding set,
\begin{equation}\label{1-RChe-6ab}
\Tg_o := \left\{ \theta_{i,j} ~|~ 1 \le i,j \le R-2 \right\}\,.
\end{equation}\medskip

\noib For $* \in \{0,\pi\}$, and $1 \le j \le R-2\,$, define
$\theta(*,q_j)$, alias $\theta_{*,j}$, to be the unique angle in the
interval $[0,\pi[$ such that
\begin{equation}\label{1-RChe-6b}
\cos\theta_{*,j} \, U_{R-1}(\cos q_j) + \sin\theta_{*,j} \,
U_{R-1}(\cos *) = 0\,.
\end{equation}
Let $\Tg_x$ denote the corresponding set,
\begin{equation}\label{1-RChe-6bb}
\Tg_x := \left\{ \theta_{*,j} ~|~ * \in \{0,\pi\}, ~1 \le j \le R-2
~ \right\}.
\end{equation} \medskip

\noib For $1 \le i \le R-2$, and $* \in \{0,\pi\}$, define
$\theta(q_i,*)$, alias $\theta_{i,*}$\,, to be the unique angle in the
interval  $[0,\pi[$ such that
\begin{equation}\label{1-RChe-6c}
\cos\theta_{i,*} \, U_{R-1}(\cos *) + \sin\theta_{i,*} U_{R-1}(\cos
q_i) = 0\,.
\end{equation}
Let $\Tg_y$ denote the corresponding set,
\begin{equation}\label{1-RChe-6cb}
\Tg_y := \left\{ \theta_{i,*} ~|~ * \in \{0,\pi\}\,, ~1 \le i \le
R-2\right\}.
\end{equation}

Observe the following relations between the above values of
$\theta$,
\begin{equation}\label{1-RChe-7c}
\theta(q_j,q_i) = \frac{\pi}{2} - \theta(q_i,q_j)\,.
\end{equation}
When $R = 2r+1$ is odd, we have
\begin{equation}\label{1-RChe-7co}
\begin{split}
\theta(q_i,q_j) = \theta(\pi-q_i,\pi-q_j) &= \theta(\pi-q_i,q_j) =
\theta(q_i,\pi-q_j)\,,\\
\theta(0,q_j) &= \theta(\pi,q_j)\,,\\
\theta(q_i,0) &= \theta(q_i,\pi)\,.
\end{split}
\end{equation}
When $R = 2r$ is even, we have
\begin{equation}\label{1-RChe-7ce}
\begin{split}
\theta(q_i,q_j) &= \theta(\pi-q_i,\pi-q_j)\,,\\
\theta(\pi-q_i,q_j) &= \pi - \theta(q_i,q_j)\,,\\
\theta(q_i,\pi-q_j) &= \pi - \theta(q_i,q_j)\,,\\
\theta(\pi,q_j) &= \pi - \theta(0,q_j)\,,\\
\theta(q_i,\pi) &= \pi - \theta(q_i,0)\,.
\end{split}
\end{equation}\bigskip

Finally, define the number $\theta_{-}$ to be,
\begin{equation}\label{1-RChe-8}
\theta_{-} := \arctan \big( \frac{1}{R}|\inf_{[-1,1]}U_{R-1}| \big)\,.
\end{equation}
We have $0 < \theta_{-} \le \pi/4\,$, with $\theta_{-} = \pi/4$ when
$R$ is even, and $\theta_{-} < \pi/4$ when $R$ is odd.\medskip

\textbf{Remark}. The pictures and numerical computations seem
to indicate that the infimum is achieved at $\cos q_1$.\medskip

\textbf{Examples.} Numerical computations give the following
approximate data when $R=8$ or $R=9$. The indication $\pi$ after the
set means that the values in the set should be multiplied by $\pi$.

$\bullet$~Special values of $\theta$ when $R=8\,$.
\begin{equation}\label{1-RChe-EX8}
\begin{split}
\Qg &= \{0.179749, 0.309108, 0.436495, 0.563505, 0.690892,
0.820251\}\,\pi ,\\
\Tg_o &= \{0.161605, 0.185335, 0.223323, 0.25, 0.276677, 0.314665,
0.338395, \\ &\hspace{1cm}0.661605, 0.685335, 0.723323, 0.75,
0.776677, 0.814665,
0.838395\}\,\pi ,\\
\Tg_x &= \{0.040363, 0.047665, 0.071705, 0.928295, 0.952335,
0.959636\}\,\pi ,\\
\Tg_y &= \{0.428295, 0.452335, 0.459636, 0.540363, 0.547665,
0.571705\}\,\pi\, .
\end{split}
\end{equation}

$\bullet$~Special values of $\theta$ when $R=9\,$.
\begin{equation}\label{1-RChe-EX9}
\begin{split}
\Qg &= \{0.159593, 0.274419, 0.387439, 0.500000, 0.612561, 0.725581,
0.840407\}\,\pi\, ,\\
\Tg_o &= \{0.145132, 0.181901, 0.217145, 0.239975, 0.260025\,,
0.282855, \\ &\hspace{1cm}0.318099, 0.354868, 0.653215, 0.707395,
0.75, 0.792605,
0.846785\}\, \pi\, ,\\
\Tg_x &= \{0.037494, 0.070922, 0.953949, 0.964777\}\, \pi\, ,\\
\Tg_y &= \{0.429078, 0.462505, 0.535223, 0.546050\}\, \pi\, .
\end{split}
\end{equation}

Up to symmetries, one can actually reduce the range of the parameter
$\theta$ to $[0,\pi/4]$ when $R$ is even, and to $[\pi/4, 3\pi/4]$
when $R$ is odd, see Subsection~\ref{SS-1RSym}. Up to this
reduction, the above values correspond to the values which appear in
the figures showing the nodal patterns for the eigenvalues
$\hat{\lambda}_{1,8}$ and $\hat{\lambda}_{1,9}$, see
Figures~\ref{Plate-1-8}, \ref{Plate-1-9-1}, and \ref{Plate-1-9-2}
at the end of the paper).

\subsection{Symmetries of the eigenfunctions associated with
$\hat{\lambda}_{1,R}$}\label{SS-1RSym}

When studying the family of eigenfunctions $\{\Phi^{\theta}\}$
associated with the eigenvalue $\hat{\lambda}_{1,R}$, it is useful
to take symmetries into account.\bigskip

\begin{properties}\label{1RSym-P1}
The following relations hold for any $(x,y) \in
[0,\pi]\times[0,\pi]$ and $\theta \in [0,\pi[\,$.
\begin{enumerate}
\item For any $ R\in \nzb$,
\begin{equation}\label{1RSym-1}
\Phi(\pi-x,\pi-y,\theta) = (-1)^{R+1}\Phi(x,y,\theta)\,.
\end{equation}
This relation implies that the nodal set $N(\Phi^{\theta})$ is
symmetrical with respect to the center $O$ of the square $\Sg$.
Furthermore,
\begin{equation}\label{1RSym-2}
\Phi(x,y,\frac{\pi}{2}-\theta) = \Phi(y,x,\theta)\,.
\end{equation}
\item When $R$ is odd, the function $\Phi$ has more symmetries,
namely,
\begin{equation}\label{1RSym-3o}
\Phi(\pi-x,y,\theta) = \Phi(x,\pi-y,\theta) = \Phi(x,y,\theta)\,.
\end{equation}
This means that the nodal set $N(\Phi^{\theta})$ is symmetrical with
respect to the lines $\{x=\pi/2\}$ and $\{y=\pi/2\}\,$.
\item When $R$ is even, we have
\begin{equation}\label{1RSym-3e}
\Phi(x,\pi-y,\theta) = \Phi(x,y,\pi-\theta) = -
\Phi(\pi-x,y,\theta)\,.
\end{equation}
\item Up to symmetries with respect to the first diagonal, or to the lines
$\{x=\pi/2\}$ and $\{y=\pi/2\}$, the nodal patterns of the family of
eigenfunctions $\{\Phi^{\theta}\}$, are those displayed by the
sub-families $\theta \in [0,\pi/4]$ when $R$ is even, and $\theta
\in [\pi/4,3\pi/4]$ when $R$ is odd.
\end{enumerate}
\end{properties}

\subsection{Critical zeroes of the eigenfunctions associated with
$\hat{\lambda}_{1,R}$}\label{SS-1RCZ}

Recall that a \emph{critical zero} of the eigenfunction
$\Phi^{\theta}$ is a point $(x,y) \in \overline{\Sg}$ such that
\begin{equation}\label{1RCZ-1}
\Phi(x,y,\theta) = \Phi_x(x,y,\theta) = \Phi_y(x,y,\theta) = 0\,.
\end{equation}

At a critical zero, the nodal set $N(\Phi^{\theta})$ consists of
several arcs (or semi-arcs when the point is on $\partial \Sg$)
which form an equi-angular system, see Properties~\ref{G-P1}. Away
from the critical zeroes, the nodal set consists of smooth embedded
arcs. To determine the possible critical zeroes of $\Phi^{\theta}$
is the key to describing the global aspect of the nodal set
$N(\Phi^{\theta})$. \medskip

We classify the critical zeroes into three (possibly empty)
categories : (i) the \emph{vertices} of the square $\Sg$; (ii) the
\emph{edge critical zeroes} located in the interior of the edges,
typically a point of the form $(0,y)$, with $y \in ]0,\pi[$; (iii)
the \emph{interior critical zeroes} of the form $(x,y) \in \Sg$.

\subsubsection{Behaviour at the vertices}\label{SSS-CZV}

Using the symmetry of $N(\Phi)$ with respect to the point $O$, see
\eqref{1RSym-1}, it suffices to consider the vertices $(0,0)$ and
$(0,\pi)$. Recalling \eqref{1R-1} and  \eqref{1RChe-2},  the Taylor expansion at $(0,0)$  of $\phi(x,y,\theta)$ is
given by,
\begin{equation}\label{1RCZV-1b}
\begin{split}
\phi(x,y,\theta) = & R \left( \cos\theta + \sin\theta \right) \, \\
& +  \frac{R(1-R^2)}{6} \,   \, \left( \cos\theta\,
y^2 + \sin\theta\,  x^2 \right)\\
& +  \, O(x^4+y^4) \,.\\
\end{split}
\end{equation}


When $R$ is odd, the behaviour is the same at the four vertices and
given by \eqref{1RCZV-1b}, due to the symmetries
\eqref{1RSym-3o}.\medskip

When $R$ is even, the Taylor expansion of $\phi(x,y,\theta)$ at
$(0,\pi)$, follows from the previous one and relation
\eqref{1RSym-3e}. In the variables $x$ and $z$ such that $y=\pi-z$,
we have,
\begin{equation}\label{1RCZV-2n}
\begin{split}
\phi(x,\pi-z,\theta) = & R \left( - \cos\theta + \sin\theta \right) \, \\
& - \frac{R(1-R^2)}{6} \,   \, \left(
-\cos\theta\,
z^2 + \sin\theta\, x^2 \right)\\
& +  \, O(x^4 + z^4) \,.\\
\end{split}
\end{equation}

\begin{figure}[!ht]
\begin{center}
\includegraphics[width=12cm]{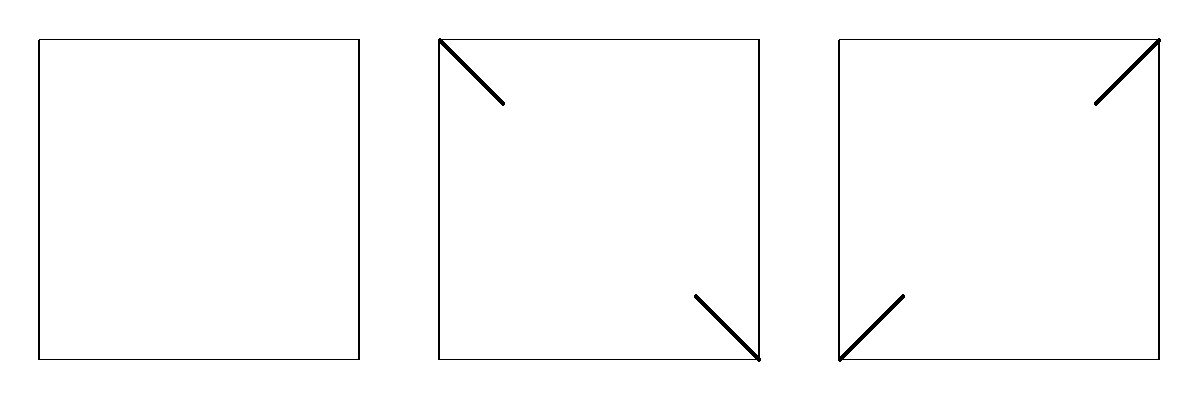}
\caption{$R$ even, nodal patterns at the vertices. From left to
right, $\theta \not = \pi/4$ and $3\pi/4$~; $\theta = \pi/4$~;
$\theta = 3\pi/4\,$} \label{1RCZV-F1e}
\end{center}
\end{figure}

\begin{figure}[!ht]
\begin{center}
\includegraphics[width=12cm]{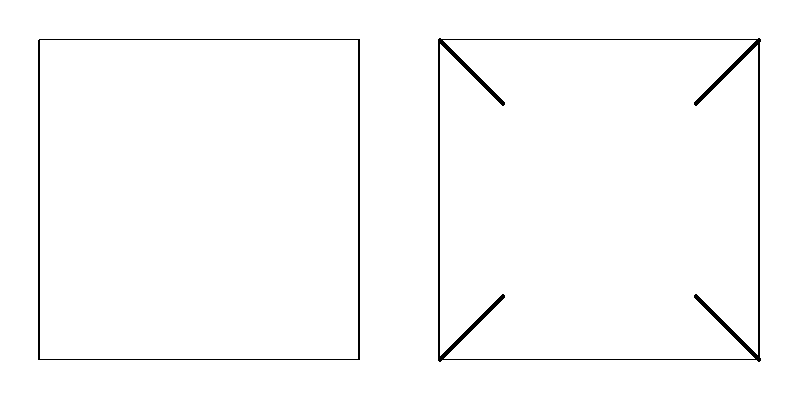}
\caption{$R$ odd, nodal patterns at the vertices. From left to
right, $\theta \not = 3\pi/4$~; $\theta = 3\pi/4\,$}\label{1RCZV-F1o}
\end{center}
\end{figure}

With the link between $\Phi$ and $\phi$ in mind, we obtain:
\begin{properties}\label{1RCZV-P1}
The vertices of the square $\Sg$ are critical zeroes for the
eigenfunction $\Phi^{\theta}$ for all $\theta$. \vspace{-3mm}
\begin{enumerate}
    \item \textbf{Case $R$ even}. The vertices $(0,\pi)$ and $(\pi,0)$ are
    non degenerate critical zeroes of $\Phi^{\theta}$ if and only if $\theta \not = \pi/4\,$.
    When $\theta = \pi/4$, they are degenerate critical zeroes of
    order $4$. The vertices $(0,0)$ and $(\pi,\pi)$ are non-degenerate
    critical zeroes of $\Phi$ if and only if $\theta \not = 3\pi/4\,$.
    When $\theta = 3\pi/4\,$, they are degenerate critical zeroes of
    order $4$. The nodal patterns at the vertices are shown in
    Figure~\ref{1RCZV-F1e}.
    \item \textbf{Case $R$ odd}. The four vertices are non-degenerate
    critical zeroes of $\Phi^{\theta}$ if and only if $\theta \not = 3\pi/4\,$.
    When $\theta = 3\pi/4\,$, they are degenerate critical zeroes of
    order $4$. The nodal patterns at the vertices are shown in
    Figure~\ref{1RCZV-F1o}.
\end{enumerate}
\end{properties}

\subsubsection{Critical zeroes, formulas}\label{SSS-CZF}

To determine the critical zeroes of the eigenfunction $\Phi$, we
 recall our notation at the beginning of the section. The first partial derivatives with respect to $x$ and $y$ are given
by,
\begin{equation}\label{1RCZF-2}
\begin{split}
\Phi_x(x,y,\theta) = & \cos x\, \sin y \; \phi(x,y,\theta)\\
&- \sin\theta \sin^2 x\, \sin y \; U_{R-1}'(\cos x)\,,\\[5pt]
\Phi_y(x,y,\theta) = & \sin x \,  \cos y \, \phi(x,y,\theta)\\
&- \cos\theta\,  \sin x\,  \sin^2 y  \;U_{R-1}'(\cos y)\,.
\end{split}
\end{equation}
The second partial derivatives are given by,
\begin{equation}\label{1RCZF-3}
\begin{split}
\Phi_{xx}(x,y,\theta) = & - \sin x\, \sin y \; \phi(x,y,\theta)\\
&- 3 \sin\theta \, \cos x\,  \sin x \, \sin y \; U_{R-1}'(\cos x)\\
& + \sin\theta \sin^3 x\, \sin y \, U_{R-1}''(\cos x)\,,\\[5pt]
\Phi_{xy}(x,y,\theta) = & \cos x\,  \cos y   \; \phi(x,y,\theta)\\
&- \cos\theta\, \cos x \, \sin^2 y \; U_{R-1}'(\cos y)\\
& - \sin\theta \, \sin^2 x\,  \cos y \; U_{R-1}'(\cos x)\,,\\[5pt]
\Phi_{yy}(x,y,\theta) = & - \sin x\,  \sin y \, \phi(x,y,\theta)\\
&- 3 \cos\theta \, \sin x \, \cos y\,  \sin y \; U_{R-1}'(\cos y)\\
& + \cos\theta \, \sin x\, \sin^3 y \;U_{R-1}''(\cos y)\,.\\
\end{split}
\end{equation}

\subsubsection{Behaviour along the edges}\label{SSS-CZE}

Recall the notation of Subsection~\ref{SS-1RChe}. Due \eqref{1RSym-1},
the symmetry with respect to the center $O$ of the square $\Sg$, it
suffices to consider the open edges $\{0\}\times]0,\pi[$ and
$]0,\pi[\times\{0\}$.\medskip

\noib \textbf{Critical zeroes on the edge} $\{0\}\times]0,\pi[\,$.
Using formulas \eqref{1RCZF-2} and \eqref{1RCZF-3}, as well as
Properties~\ref{1RChe-P1}, we infer that the point $(0,y)$ is a
critical zero for $\Phi^{\theta}$ if and only if,
\begin{equation}\label{1RCZE-1y}
\cos\theta \, U_{R-1}(\cos y) + R \sin\theta = 0\,,
\end{equation}
with second derivatives at $(0,y)$, $\Phi_{xx} = \Phi_{yy} = 0\,$, and
$ \Phi_{xy} = - \cos\theta\,  \sin^2 y \; U_{R-1}'(\cos y)\,.$

The point $(0,y)$ is a non degenerate critical zero, unless $y = q_j
\in \Qg$ for some $j, 1 \le j \le R-2$. This can only occur when
$\theta = \theta(0,q_j)$. In this case, the third derivative
$\Phi_{xy^2}$ at the degenerate critical zero $(0,q_j)$ is equal to
$$ \cos\left(\theta(0,q_j)\right)\, \sin^3 q_j\; U_{R-1}''(\cos q_j) \not = 0\,,$$
and the critical zero has order $3$.
\medskip

\noib \textbf{Critical zeroes on the edge} $]0,\pi[\times\{0\}$.
Similarly, the point $(x,0)$ is a critical zero for $\Phi^{\theta}$
if and only if,
\begin{equation}\label{1RCZE-1x}
R\,\cos\theta  + \sin\theta\, U_{R-1}(\cos x) = 0\,,
\end{equation}
with second derivatives at $(x,0)$, $\Phi_{xx} = \Phi_{yy} = 0\,$, and
$$ \Phi_{xy} = - \sin\theta \, \sin^2 x \; U_{R-1}'(\cos x)\,.$$

The point $(x,0)$ is a non degenerate critical zero unless $x = q_i
\in \Qg$ for some $i, 1 \le i \le R-2\,$. This can only occur when
$\theta = \theta(q_i,0)$. In this case, the third derivative
$\Phi_{x^2y}$ at the degenerate critical zero $(q_i,0)$ is equal to
$$
\sin\left(\theta(q_i,0)\right) \, \sin^3 q_i \; U_{R-1}''(\cos q_i)
\not = 0\,,
$$
and the critical zero has order $3$. \medskip

\textbf{Remark}. At an edge critical zero which is non
degenerate, an arc from the nodal set hits the edge orthogonally. At
a degenerate edge critical zero, two arcs from the nodal set hit the
edge with equal angle $\pi/3$. See  Figures~\ref{1RQ-F2e} and
\ref{1RQ-F2o}.
\medskip

The following properties summarize the analysis of the above
equations.

\begin{properties}\label{1RCZE-P1}
The critical zeroes on the open edges, if any, appear in pairs of
points which are symmetrical with respect to the center $O$ of the
square $\Sg$.\\[4pt]
\textbf{Case $R$ even}.\vspace{-3mm}
\begin{enumerate}
    \item For $\theta \in [0,\pi/4[ \cup
    ]3\pi/4,\pi[$, there are critical zeroes on the vertical edges
    $\{0,\pi\}\times]0,\pi[$, and no critical zero on the horizontal
    edges $]0,\pi[\times\{0,\pi\}$.
    \item For $\theta \in ]\pi/4,3\pi/4[$, there are critical zeroes on the
    horizontal edges $]0,\pi[\times\{0,\pi\}$, and no critical zero on the
    vertical edges $\{0,\pi\}\times]0,\pi[\,$.
    \item The number of critical zeroes depends on $\theta$, more precisely on
    the number of solutions of \eqref{1RCZE-1y} or \eqref{1RCZE-1x}.
    \item When $\theta = \pi/4$ or $3\pi/4\,$, the only boundary
    critical zeroes are vertices, see Properties~\ref{1RCZV-P1}.
\end{enumerate}\vspace{-3mm}
\textbf{Case $R$ odd}.\vspace{-3mm}
\begin{enumerate}
    \item Recall the value $0 < \theta_{-} < \pi/4$ defined in Subsection~\ref{SS-1RChe}.
    For $\theta \in [0,\theta_{-}] \cup ]3\pi/4,\pi[\,$, there are critical zeroes on the
    vertical edges $\{0,\pi\}\times]0,\pi[\,$, and no critical zero on the horizontal edges
    $]0,\pi[\times\{0,\pi\}$.
    \item For $\theta \in [\pi/2-\theta_{-},3\pi/4[\,$, there are critical zeroes on the
    horizontal edges $]0,\pi[\times\{0,\pi\}$, and no critical zero on the vertical
    edges $\{0,\pi\}\times]0,\pi[\,$.
    \item For $\theta \in ]\theta_{-}, \pi/2-\theta_{-}[\,$, there is no critical zero on the
    open edges.
    \item  The number of critical zeroes depends on $\theta$, more precisely on
    the number of solutions of \eqref{1RCZE-1y} or \eqref{1RCZE-1x}.
    \item The critical zeroes have order
    at most $3$. Degenerate critical zeroes can only occur
    for finitely many values of $\theta$ and $x$ or $y\,$.
    \item When $\theta = 3\pi/4\,$, the only boundary
    critical zeroes are the vertices, see Properties~\ref{1RCZV-P1}.
\end{enumerate}\vspace{-3mm}
In both cases, the edge critical zeroes are non degenerate unless
they occur on a horizontal edge for some $x=q_i \in \Qg$, \resp on a
vertical edge for some $y = q_j \in \Qg$, in which case $\theta$
must be equal to $\theta(q_i,0)$, \resp to $\theta(0,q_j)$.
Degenerate critical zeroes have order $3$. If $\theta \not = 0$ or
$\pi/2\,$, the points $(*,p_j)$ and $(p_j,*)$ with $1 \le j \le R-1$
and $* = 0$ or $\pi$ are not critical zeroes of the eigenfunction
$\Phi^{\theta}$.
\end{properties}

\textbf{Remark}. A more detailed description of the
localization of the edge critical zeroes is given in
Subsection~\ref{SS-1RQ}.

\subsubsection{Interior critical zeroes}\label{SSS-CZI}

Recall the notations of Subsection~\ref{SS-1RChe}. The following
properties follow from \eqref{1RCZF-2}--\eqref{1RCZF-3}.

\begin{properties}\label{1RCZI-P1} Let $(x,y) \in \Sg$ be an interior
point. \vspace{-3mm}
\begin{enumerate}
    \item The functions $\Phi$ and $\Phi_x$ vanish at $(x,y)$ if and only if
    \begin{equation}\label{1RCZI-1x}
    \begin{split}
    & \cos\theta \, U_{R-1}(\cos y) + \sin\theta \, U_{R-1}(\cos x) = 0\,,
    \text{~and}\\
    & U_{R-1}'(\cos x) = 0\,.\\
    \end{split}
    \end{equation}
    This happens in particular at regular points of the nodal set with
    a horizontal tangent.
    \item The functions $\Phi$ and $\Phi_y$ vanish at $(x,y)$ if and only if
    \begin{equation}\label{1RCZI-1y}
    \begin{split}
    & \cos\theta \, U_{R-1}(\cos y) + \sin\theta  \, U_{R-1}(\cos x) = 0\,,
    \text{~and}\\
    & U_{R-1}'(\cos y) = 0\,.\\
    \end{split}
    \end{equation}
    This happens in particular at regular points of the nodal set with
    a vertical tangent.
    \item The point $(x,y)$ is an interior critical zero of
    $\Phi$, if and only if
    \begin{equation}\label{1RCZI-1}
    \begin{split}
    & \cos\theta \, U_{R-1}(\cos y) + \sin\theta \,  U_{R-1}(\cos x)
    = 0\,, \text{~and~}\\
    & U_{R-1}'(\cos x) = 0 \text{~and~} U_{R-1}'(\cos y) = 0\,.
    \end{split}
    \end{equation}
    The only possible interior critical zeroes for the family of
    eigenfunctions $\{\Phi^{\theta}\}$ are the points $(q_i,q_j), 1 \le
    i,j \le R-2\,$. The point $(q_i,q_j)$ can only occur as a critical
    zero of the eigenfunction $\Phi^{\theta(q_i,q_j)}\,$.
    When $\theta$ is not one of the values $\theta(q_i,q_j)$, $1\le i,j\le R-2\,$,
    the eigenfunction
    $\Phi^{\theta}$ does not have any interior critical zero.
    \item When $(x,y)$ is an interior critical zero of $\Phi$, the
    Hessian of $\Phi$ at $(x,y)$ is given by,
    \begin{equation*}
    \sin x \, \sin y \,\begin{pmatrix}
    \sin\theta \,\sin^2 x\; U_{R-1}''(\cos x) & 0 \\
    0 & \cos\theta\, \sin^2 y\; U_{R-1}''(\cos y) \\
    \end{pmatrix},
    \end{equation*}
    so that the interior critical zeroes, if any, are always non
    degenerate.
    \item The lattice points $\Lg = \{(\frac{i\pi}{R},\frac{j\pi}{R})\,, ~1
    \le i,j \le R-1\}$ (see Subsection~\ref{SS-GNot2}) are common zeroes to
    all the eigenfunctions $\Phi^{\theta}$ when $\theta \in [0,\pi[\,$. They are
    not interior critical zeroes.
\end{enumerate}
\end{properties}

\subsection{$Q$-nodal patterns of eigenfunctions associated with
$\hat{\lambda}_{1,R}$}\label{SS-1RQ}

The purpose of this section is to list all the possible patterns of
the nodal set $N(\Phi^{\theta})$ inside the $Q$-squares $Q_{i,j}$,
$0 \le i,j \le R-1$, see Subsection~\ref{SS-GNot2}.

The following properties are derived from the previous sections and
from Properties~\ref{G-P2}.\vspace{-3mm}
\begin{enumerate}
    \item The nodal set $N(\Phi^{\theta})$ is contained in $\Wg(\pm) \cup \Lg
    \cup \partial \Sg$.
    \item If a white square $Q_{i,j}$ does not touch the boundary
    $\partial \Sg$, the nodal set $N(\Phi^{\theta})$ cannot cross nor hit the boundary
    of $Q_{i,j}$, except at the vertices which belong to $\Lg$.
    \item If a white square $Q_{i,j}$ touches the boundary
    $\partial \Sg$, the nodal set $N(\Phi^{\theta})$ cannot intersect the boundary
    of $Q_{i,j}$, except at the vertices which belong to $\Lg$, or at
    an edge contained in $\partial \Sg$.
    \item Since the points in $\Lg$ are not critical zeroes, the
    nodal set consists of a single regular arc at such a point.
    \item Inside a white square $Q_{i,j}$, the nodal set $N(\Phi^{\theta})$
    can have at most one self intersection at $(q_i,q_j)$ if this
    point is a critical zero for $\Phi^{\theta}$. In this case, the critical
    zero is non degenerate, and the nodal set at $(q_i,q_j)$ consists
    locally of two regular arcs meeting orthogonally.
    \item Inside a square $Q_{i,j}$, the nodal set cannot stop at a
    point, and consists of at most finitely many arcs.
    \item The nodal set $N(\Phi^{\theta})$ cannot contain a closed curve
    contained in the closure of $Q_{i,j}$ (energy reasons).
\end{enumerate}

\begin{figure}[!ht]
\begin{center}
\includegraphics[width=12cm]{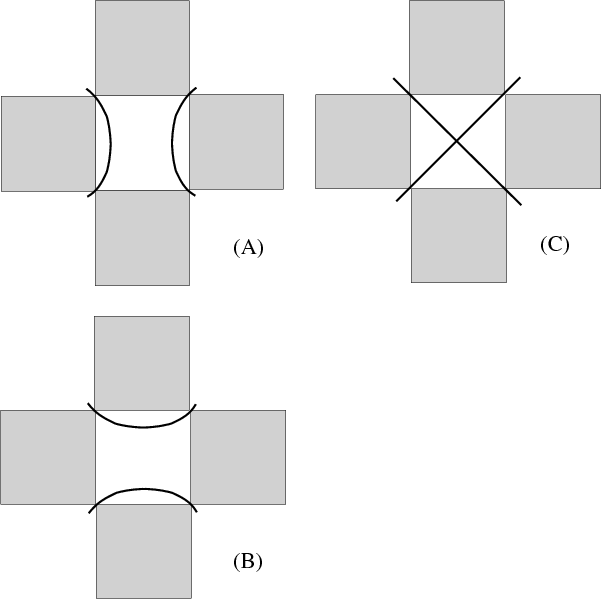}
\caption{Nodal pattern in an inner $Q$-square} \label{1RQ-F1}
\end{center}
\end{figure}

$\bullet$~\textbf{Inner $Q$-square}. Figure~\ref{1RQ-F1} shows all
the possible nodal patterns inside a square $Q_{i,j}$ which does not
touch the boundary. The patterns~A and B occur when the
eigenfunction $\Phi^{\theta}$ does not have any interior critical
zero inside the $Q$-square. Pattern~C occurs when $\Phi^{\theta}$
admits $(q_i,q_j)$ as interior critical zero (necessarily unique and
non degenerate), in which case $\theta$ must be equal to
$\theta(q_i,q_j)$. The properties recalled above show that there are
no other possible nodal patterns. \medskip

$\bullet$~\textbf{Boundary $Q$-square, $R$ even}.  As stated in
Properties~\ref{1RSym-P1}~(iv), it suffices to consider the case
$\theta \in [0,\pi/4]$. The only boundary critical zeroes of
$\Phi^{\theta}$ are the vertices $(0,\pi)$ and $(\pi,0)$, this case
occurs if and only if $\theta = \pi/4\,$, and points on the vertical
edges $\{0,\pi\}\times]0,\pi[$ if and only if $0\le \theta <
\pi/4\,$. These points come in pairs of symmetric points with
respect to the center $O$ of the square $\Sg$. It suffices to
describe the points located on the edge $\{0\}\times]0,\pi[\,$,
 \ie the points $(0,y)$ satisfying equation \eqref{1RCZE-1y},
$$U_{R-1}(\cos y) + R\tan\theta =0\,, \text{~for some~} \theta\,,
~~0\le \theta < \pi/4\,.$$
When $\theta = 0$, this equation provides exactly $(R-1)$ non
degenerate critical zeroes, the points $p_j, 1\le j \le R-1\,$. When
$0 < \theta < \pi/4$ the equation has at least one solution located
in the interval $]p_{R-1},\pi[$. This is the sole solution when
$\theta$ is close enough to $\pi/4\,$, and it corresponds to a non
degenerate critical zero. The other solutions, if any, are located
in the intervals $]p_j,p_{j+1}[$, with $j$ odd. Each such interval
contains at most two solutions, which correspond to non degenerate
critical zeroes. When an interval contains only one point, this is a
double solution, and it corresponds to a degenerate critical zero
$q_j$. This can only occur for the special value
$\theta(0,q_j)\,$.\medskip

\begin{figure}[!ht]
\begin{center}
\includegraphics[width=12cm]{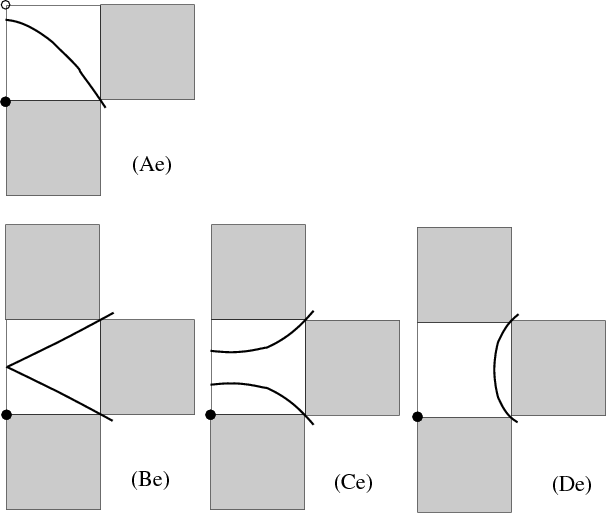}
\caption{Local nodal patterns at the boundary, $R$
even}\label{1RQ-F2e}
\end{center}
\end{figure}

Figure~\ref{1RQ-F2e} gives all the possible nodal patterns of the
eigenfunction $\Phi^{\theta}$ in a square $Q_{0,j}$ which touches
the edge $\{0\}\times]0,\pi[$. The base point of the square
$(0,p_j)$ is the black dot in the figures. The other dot is the
vertex $(0,\pi)$.\medskip

Figure~(Ae) shows the nodal pattern in the square $Q_{0,R-1}$ which
touches the vertex $(0,\pi)$, and contains the persistent edge non
degenerate critical zero. Figure~(Be) shows the nodal pattern when
there is a degenerate edge critical zero $(0,q_j)$ ; it is always of
order three, with two arcs hitting the boundary with equal angles
$\pi/3$. Figure~(Ce) shows the nodal pattern when there are two non
degenerate critical zeroes in the interval
$\{0\}\times]p_j,p_{j+1}[\,$. There are two arcs hitting the boundary
orthogonally. Figure~(De) shows the nodal pattern when the interval
$\{0\}\times]p_j,p_{j+1}[$ contains no critical zero. The properties
recalled above show that there are no other possible nodal patterns.
\bigskip

\begin{figure}[!ht]
\begin{center}
\includegraphics[width=12cm]{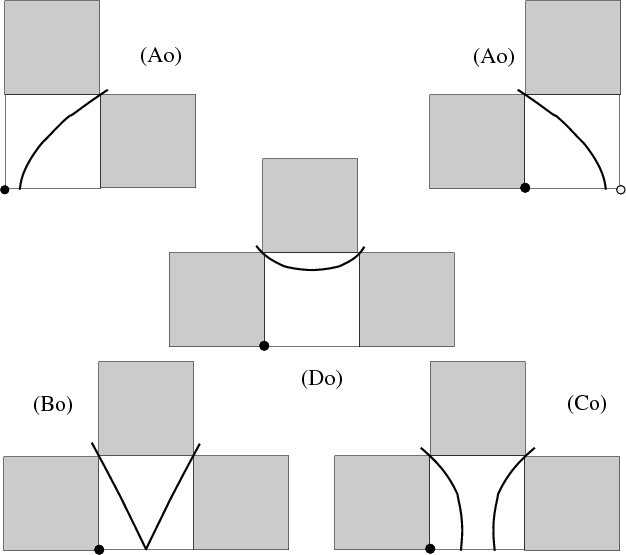}
\caption{Local nodal patterns at the boundary, $R$
odd}\label{1RQ-F2o}
\end{center}
\end{figure}

$\bullet$~\textbf{Boundary $Q$-square, $R$ odd}.  The description of
the critical zeroes of $\Phi^{\theta}$ on the boundary $\partial
\Sg$ in the case $R$ odd is similar to the case $R$ even, with some
changes. As stated in Properties~\ref{1RSym-P1}~(iv), up to
symmetries, we can restrict ourselves to $\theta \in [\pi/4,
3\pi/4]\,$. Recall the value $0 < \theta_{-}<\pi/4$ defined in
Subsection~\ref{SS-1RChe}. The vertices are critical zeroes, and they
are non degenerate unless $\theta = 3\pi/4$. For $\theta \in
[\pi/4,\pi/2-\theta_{-}[\,$, there is no critical zero on the edges.
For $\theta \in ]\pi/2-\theta_{-},3\pi/4[\,$, there are critical
zeroes on the horizontal edges, and none on the vertical edges.
Since the nodal sets are symmetrical with respect to $\{y=\pi/2\}$,
it suffices to describe the critical zeroes on the horizontal edge
$]0,\pi[\times\{0\}$. For $\theta = \pi/2-\theta_{-}\,$, there are at
least two order $3$ critical zeroes (except when $R=3$ in which case
there is only one). As a matter of fact, it seems that there are
exactly two critical zeroes for $R \ge 5$ because the local extrema
of $U_{R-1}$ decrease in absolute value on $[-1,0]$. For $\theta \in
]\pi/2-\theta_{-},\pi/2]\,$, the number of critical zeroes depends on
the number of solutions of equation \eqref{1RCZE-1x},
$$R\,\cot\theta + U_{R-1}(\cos x) = 0\,,$$
with $0$ or $2$ solutions in each interval
$]p_i,p_{i+1}[\times\{0\}\,$, or a degenerate critical zero at some
$(q_i,0)$, when $\theta = \theta(q_i,0)$. For $\theta \in
]\pi/2,3\pi/4[\,$, there are two non degenerate critical zeroes, one
in each interval $]0,p_1[\times\{0\}$ and
$]p_{R-1},\pi[\times\{0\}$, near the vertices. For $\theta =
3\pi/4\,$, there is no critical zero on the open edge, and only
critical zeroes of order $4$ at the vertices. Using the properties
listed at the beginning of Subsection~\ref{SS-1RQ}, one can show that
Figure~\ref{1RQ-F2o} contains all the possible nodal patterns in a
$Q$-square touching the edge $]0,\pi[\times\{0\}\,$.

\subsection{Nodal sets of the functions $Z_{\pm}$ associated with
$\hat{\lambda}_{1,R}$}\label{SS-1RZ}

Recall the  notation in Subsection~\ref{SS-GNot2}, \eqref{1RNot-9},  $Z_{+}= \Phi^{\frac{\pi}{4}} \text{~and~}
Z_{-}= \Phi^{\frac{3\pi}{4}}\,.$ The purpose of this section is to prove the following proposition.

\begin{proposition}\label{1RZ-P0}Nodal sets of the eigenfunctions $Z_{\pm}$
associated with the eigenvalue $\hat{\lambda}_{1,R}$. \vspace{-3mm}
\begin{enumerate}
\item \textbf{Case $R$ even, $R=2r$}. The symmetries with respect to the lines
$\{x=\pi/2\}$ and $\{y=\pi/2\}$ send the nodal set $N(Z_{+})$ to the
nodal set $N(Z_{-})$. Both nodal sets $N(Z_{+})$ and $N(Z_{-})$ are
invariant under the symmetry with respect to the center $O$ of the
square. The nodal set $N(Z_{+})$ consists of the boundary $\partial
\Sg$, the anti-diagonal, and a collection of $(r-1)$ simple closed
curves $\gamma_i$. The curve $\gamma_i$ winds around $\partial
\Sg_i$, passing through the points in $\Lg \cap \partial \Sg_i$, and
crosses the anti-diagonal $\Dg_{-}$ orthogonally. The curves
$\gamma_i$ do not intersect each other.
\item \textbf{Case $R$ odd, $R=2r+1$}. The nodal sets $N(Z_{\pm})$ are
invariant under the symmetries with respect to the lines
$\{x=\pi/2\}$ and $\{y=\pi/2\}$. The nodal set $N(Z_{+})$ consists
of the boundary $\partial \Sg$, and a collection of $r$ simple
closed curves $\alpha_i$. The curve $\alpha_i$ winds around
$\partial \Sg_i$, passing through the points in $\Lg \cap \partial
\Sg_i$. The curves $\alpha_i$ do not intersect each other. The nodal
set $N(Z_{-})$ consists of the boundary $\partial \Sg$, the two
diagonals $\Dg_{\pm}$, and a collection of $(r-1)$ simple closed curves
$\beta_i$. The curve $\beta_i$ winds around $\partial \Sg_i$,
passing through the points in $\Lg \cap \partial \Sg_i\,$, and crosses
the diagonals $\Dg_{\pm}$ orthogonally. The curves $\beta_i$ do not
intersect each other.
\end{enumerate}
\end{proposition}

The proof of this proposition relies on three lemmas which we prove
below. Lemma~\ref{1RZ-L1} determines precisely the interior critical
zeroes of the nodal sets $N(Z_{\pm})$. Lemma~\ref{1RZ-L2a} and
\ref{1RZ-L2b} are ``separation'' lemmas.\medskip

\begin{lemma}\label{1RZ-L1}
When $R$ is even, the diagonal, \resp the anti-diagonal, is
contained in the nodal set $N(Z_{-})$, \resp in the nodal set
$N(Z_{+})$. When $R$ is odd, the diagonal and anti-diagonal meet the
nodal set $N(Z_{+})$ at finitely many points. They are both
contained in the nodal set $N(Z_{-})$. The critical zeroes of the
functions $Z_{\pm}$ are as follows. \vspace{-4mm}
\begin{enumerate}
    \item \textbf{Case $R$ even}, $R=2r\,$. The interior critical zeroes of the
    function $Z_{+}$ are exactly the $(R-2)$ points, $(q_i,\pi-q_i)$,
    for $1 \le i \le R-2\,$, located on the anti-diagonal.
    The interior critical zeroes of the function $Z_{-}$ are precisely the
    $(R-2)$ points, $(q_i,q_i)$, for $1 \le i \le R-2$, located on the
    diagonal.
    \item \textbf{Case $R$ odd}, $R=2r+1\,$. The function $Z_{+}$ has no
    interior critical zero. The interior critical zeroes of the
    function $Z_{-}$ are precisely the $(2R-5)$ points,
    $(q_i,q_i)$, $(q_{i}, \pi-q_{i})$, for $1 \le i \le R-2\,$, located on the
    diagonal and anti-diagonal.
\end{enumerate}
\end{lemma}

\textbf{Remark}. Note that Lemma~\ref{1RZ-L1},
Properties~\ref{1RCZE-P1} and Properties~\ref{1RCZV-P1} provide a
complete description of the critical zeroes of the functions
$Z_{\pm}$.\medskip

 \textbf{Proof}. The first assertions are clear. We
concentrate on the determination of the interior critical zeroes.
Let $\epsilon = \pm 1$. The point $(x,y) \in \Sg$ is a critical zero
of the function $Z_{\pm}$ if and only if $(x,y)$ is a common solution to
the following three equations.
\begin{equation}\label{1RZ-1}
\sin x \, \sin(Ry) + \epsilon\, \sin(Rx) \, \sin y = 0\,.
\end{equation}
\begin{equation}\label{1RZ-1x}
\cos x \, \sin(Ry) + \epsilon R \, \cos(Rx) \, \sin y = 0\,.
\end{equation}
\begin{equation}\label{1RZ-1y}
R\, \sin x \, \cos(Ry) + \epsilon\, \sin(Rx) \, \cos y = 0\,.
\end{equation}
Since $(x,y)$ is an interior critical zero,\eqref{1RZ-1} and
\eqref{1RZ-1x} imply equation \eqref{1RZ-2x} below ; \eqref{1RZ-1}
and \eqref{1RZ-1y} imply equation \eqref{1RZ-2y} below.
\begin{equation}\label{1RZ-2x}
\cos x \, \sin(Rx) - R\, \sin x  \, \cos(Rx) = 0\,.
\end{equation}
\begin{equation}\label{1RZ-2y}
\cos y \, \sin(Ry) - R\, \sin y \, \cos(Ry) = 0\,.
\end{equation}

Substituting $\sin(Rx)$ and $\sin(Ry)$ in \eqref{1RZ-1} using
\eqref{1RZ-2x} and \eqref{1RZ-2y}, we obtain the equation
\begin{equation}\label{1RZ-3}
\cos x \, \cos(Ry) + \epsilon\, \cos(Rx) \, \cos y = 0\,.
\end{equation}

Adding and substracting \eqref{1RZ-1} to/from \eqref{1RZ-3}, we
obtain that an interior critical zero $(x,y)$ of $Z_{\pm}$ satisfies the
system,
\begin{equation}\label{1RZ-4}
\begin{split}
& \cos(x-Ry) + \epsilon\, \cos(Rx-y) = 0\,,\\
& \cos(x+Ry) + \epsilon\, \cos(Rx+y) = 0\,.\\
\end{split}
\end{equation}

\noib Case $\epsilon = 1$.  Modulo $2\pi$, the system
\eqref{1RZ-4} is equivalent to
\begin{equation}\label{1RZ-4a}
\begin{aligned}
(a1)~ x-Ry=\pi-Rx+y ~[2\pi] \text{~or~} &(b1)~ x-Ry=\pi+Rx-y
~[2\pi]\,,\\
\text{and} \\
(a2)~ x+Ry=\pi-Rx-y ~[2\pi] \text{~or~} &(b2)~ x+Ry=\pi+Rx+y
~[2\pi]\,.\\
\end{aligned}
\end{equation}

We have to consider four cases.

$\checkmark$~{(a1) \& (a2)}~ These conditions imply that $Rx=-x +
k\pi$ for some integer $k$. Using \eqref{1RZ-2x}, we find that
$$(-1)^{k+1}(1+R)\, \sin x \,\cos x = 0\,.$$
This implies that $x=\pi/2\,$. Similarly, using \eqref{1RZ-2y}, we
find that $y = \pi/2$. Since $(x,y)$ is a critical zero, using
Subsection~\ref{SS-1RChe}, this can only occur when $R$ is odd. On the
other hand, using \eqref{1RZ-1}, we find that $\sin(R\pi/2)=0$ which
implies that $R$ is even. The conditions (a1) and (a2) cannot occur
simultaneously.\medskip

$\checkmark$~{(a1) \& (b2)}~ These conditions imply that $x=y\,$.
Using \eqref{1RZ-1}, we find that $\sin(Rx)=0$ and hence, by
\eqref{1RZ-2x}, $\cos(Rx)=0$. The conditions (a1)  and (b2)
cannot occur simultaneously.\medskip

$\checkmark$~{(b1) \& (a2)}~ These conditions imply that $y=\pi-x\,$.
Using \eqref{1RZ-1}, we see that
$$\big( (-1)^{R+1} + 1\big) \, \sin x\, \sin(Rx) =0\,.$$ If $R$ were odd,
we would have a contradiction with \eqref{1RZ-2x}. This case can
only occur when $R$ is even. \medskip

$\checkmark$~{(b1) \& (b2)}~ These conditions imply that $Rx=x +
k\pi$ for some integer $k\,$. By \eqref{1RZ-2x}, this implies that
$x=\pi/2$. Similarly, we find that $y=\pi/2$. As above, this implies
that $R$ is odd. On the other-hand, \eqref{1RZ-1}, implies that
$\sin(R\pi/2)$ which implies that $R$ is even. The conditions (b1)
and (b2) cannot occur simultaneously.\medskip

We conclude that the function $Z_{+}$ has no interior critical zero
when $R$ is odd, and that its only critical zeroes are the points
$(q_i,\pi-q_i)$, for $1\le i \le R-2$ when $R$ is even.\bigskip

\noib Case $\epsilon = -1$. The system \eqref{1RZ-4} is equivalent
to
\begin{equation}\label{1RZ-4b}
\begin{aligned}
(a1)~ x-Ry=Rx-y ~[2\pi] \text{~or~} &(b1)~ x-Ry=-Rx+y ~[2\pi]\,,\\
\text{and} \\
(a2)~ x+Ry=Rx+y ~[2\pi] \text{~or~} &(b2)~ x+Ry=-Rx-y ~[2\pi]\,.\\
\end{aligned}
\end{equation}

We have to consider four cases.\\

$\checkmark$~{(a1) \& (a2)}~ These conditions imply that $Rx=x +
k\pi$ for some integer $k$. Using \eqref{1RZ-2x}, we find that
$$(-1)^{k+1}(1-R)\, \sin x \,\cos x = 0\,.$$
This implies that $x=\pi/2\,$. Similarly, using \eqref{1RZ-2y}, we
find that $y = \pi/2\,$. Since $(x,y)$ is a critical zero, using
Subsection~\ref{SS-1RChe}, this case can only occur when $R$ is
odd.\medskip

$\checkmark$~{(a1) \& (b2)}~ These conditions imply that $\pi-x=y\,$.
Using \eqref{1RZ-1}, we find that $$\big((-1)^{R+1}-1\big)\,
\sin(Rx)=0\,.$$
Since $(x,y)$ is a critical zero, $\sin(Rx) \not = 0$
and this case can only occur when $R$ is odd.\medskip

$\checkmark$~{(b1) \& (a2)}~ These conditions imply that $y=x\,$. This
case occurs for both $R$ even and $R$ odd. \medskip

$\checkmark$~{(b1) \& (b2)}~ These conditions imply that $Rx=-x +
k\pi$ for some integer $k$. By \eqref{1RZ-2x}, this implies that
$x=\pi/2\,$. Similarly, we find that $y=\pi/2\,$. This case can only
occur when $R$ is odd.\medskip

We conclude that the only critical zeroes of the function $Z_{-}$
are the points $(q_i,q_i)$, for $1\le i \le R-2$ when $R$ is even,
and are the points, $(q_i,q_i)$, $(q_i,\pi-q_i)$, for $1\le i \le
R-2$ when $R$ is odd.\hfill $\square$ \bigskip

Recall that the function $Z_{+}(x,y)$ satisfies the relations
$$Z_{+}(y,x) = Z_{+}(x,y) \text{~~~and~~~} Z_{+}(\pi-x,\pi-y) =
(-1)^{R+1}\,Z_{+}(x,y)$$ which imply that the nodal set $N(Z_{+})$
is invariant under the symmetry with respect to the diagonal
$\Dg_{+}$\,, and under the symmetry with respect to the centre $O$
of the square $\Sg$. Consider the subsets
\begin{equation}\label{1RZ-5}
\begin{split}
\Fg_1 & := \left\{ \Sg \cap \{x > y\} \cap \{x+y < \pi \right\},\\
\Fg_2 & := \left\{ \Sg \cap \{x > y\} \cap \{x+y > \pi \right\},\\
\Fg_3 & := \left\{ \Sg \cap \{x < y\} \cap \{x+y > \pi \right\},\\
\Fg_4 & := \left\{ \Sg \cap \{x < y\} \cap \{x+y < \pi \right\}.\\
\end{split}
\end{equation}

Due to the symmetries mentioned above, it suffices to understand the
nodal set into one of theses domains. Since $Z_{+}$ corresponds to
the value $\theta = \pi/4\,$, the diagonal $\Dg_{+}$ is covered by
grey $Q$-squares. When $R$ is even, the anti-diagonal $\Dg_{-}$ is
covered by white squares which either contain a unique critical zero
of $Z_{+}$, or have as vertex one of the vertices $(0,\pi)$ or
$(\pi,0)$. In either situations, the structure of $N(Z_{+})$ inside
these diagonal white $Q$-squares is known, see Section~\ref{SS-1RQ}.
When $R$ is odd, both diagonals $\Dg_{+}$ and $\Dg_{-}$ are covered
by grey squares, and the white $Q$-squares meeting a given $\Fg_i$
are actually contained in $\Fg_i$. In summary, it suffices to
understand the nodal pattern of $Z_{+}$ inside the white squares
contained into the $\Fg_i$, and it suffices to look at the case
$i=1$, and use the symmetries. \medskip

{\bf Claim}. \emph{In each white square $Q_{i,j} \subset \Wg(+)
\cap \Fg_1$, the horizontal segment $]p_i,p_{i+1}[\times\{m_j\}$
does not meet the nodal set $N(Z_{+})$. More precisely, for $R=2r$
and $j\le r-1\,$,
$$(-1)^{j}Z_{+}(x,m_j) > 0 \text{~on the interval~} ]p_{j+1},\pi-p_{j+1}[\,.$$}

 \textbf{Proof}. Since $Q_{i,j} \in \Wg(+)$, we must have
$i+j$ odd \ie  $(-1)^{i+j} = - 1\,$. Since $Q_{i,j} \subset
\Fg_1$, we must have the inequalities $j\le i-1$ and $i+j \le R-2\,$.
Up to the positive factor $1/\sqrt{2}\,$, we have, for any $x \in
]p_i,p_{i+1}[\,$,
\begin{equation}\label{1RZ-6a}
\begin{split}
Z_{+}(x,m_j) & = \sin x\, \sin(Rm_j) + \sin(Rx)\,\sin m_j \\
& = (-1)^{j}\big( \sin x + (-1)^j \sin m_j\, \sin(Rx) \big) \\
& = (-1)^{j}\big( \sin x - \sin m_j\,|\sin(Rx)| \big)\,,
\end{split}
\end{equation}
where the last equality follows from the equalities $\sin(Rx) =
(-1)^i|\sin(Rx)|$ on the interval $]p_i,p_{i+1}[$, and $(-1)^{i+j} =
-1$\,. On the other-hand, the inequalities $j\le i-1$ and $i+j \le
R-2$ imply that $m_j < p_i < p_{i+1} < \pi - m_j\,$, so that $\sin x
> \sin m_j $ on $]p_i,p_{i+1}[\,$. This proves that $Z_{+}(x,m_j) \not
= 0$ on $]p_i,p_{i+1}[\,$,  hence the claim. \hfill $\square$.

Taking into account the preceding discussion, we have obtained the
following lemma.

\begin{lemma}\label{1RZ-L2a}
The horizontal segments (medians) through the points $(m_i,m_j)$
which are contained in the white squares $Q_{i,j} \subset \Wg(+)
\cap \Fg_1$ or $\Wg(+) \cap \Fg_3$ do not meet the nodal set
$N(Z_{+})$. The vertical segments (medians) through $(m_i,m_j)$
which are contained in the white squares $Q_{i,j} \subset \Wg(+)\cap
\Fg_2$ or $\Wg(+)\cap \Fg_4$ do not meet the nodal set $N(Z_{+})$.
\end{lemma}

We have a similar lemma for the eigenfunction $Z_{-}$.

\begin{lemma}\label{1RZ-L2b}
The horizontal segments (medians) through $(m_i,m_j)$ which are
contained in the white squares $Q_{i,j} \subset \Wg(+) \cap \Fg_1$
or $\Wg(-) \cap \Fg_3$ do not meet the nodal set $N(Z_{-})$. The
vertical segments (medians) through $(m_i,m_j)$ which are contained
in the white squares $Q_{i,j} \subset \Wg(+)\cap \Fg_2$ or
$\Wg(+)\cap \Fg_4$ do not meet the nodal set $N(Z_{-})$.
\end{lemma}

 \textbf{Proof}. We sketch the proof of the lemma. Since
$N(Z_{+})$ and $N(Z_{-})$ are symmetrical to each other with respect
to $\{x=\pi/2\}$ when $R$ is even, it suffices to study $N(Z_{-})$
for $R$ odd. The nodal set is contained in $\Wg(-)$. The diagonal
and the anti-diagonal are covered by white $Q$-squares, and in these
squares the nodal pattern is known since they either contain a
critical zero or touch a vertex of the square $\Sg$. As above, we
look at the white squares inside $\Fg_1$. The indices of these
squares satisfy
$$(-1)^{i+j} = 1\,, ~~ j\le i-1\,, ~~ i+j\le R-2\,.$$
As in the previous proof, we can write,
\begin{equation}\label{1RZ-6b}
\begin{split}
Z_{-}(x,m_j) & = \sin x \, \sin(Rm_j) - \sin(Rx)\,\sin m_j \\
& = (-1)^{j}\big( \sin x  + (-1)^j \sin m_j \, \sin(Rx) \big) \\
& = (-1)^{j}\big( \sin x - \sin m_j|\sin(Rx)| \big)\,,
\end{split}
\end{equation}
because $\sin(Rx) = (-1)^{i}|\sin(Rx)|$ in the interval
$]p_i,p_{i+1}[\,$, and $(-1)^{i+j}=1\,$.\\
 The same argument as above
gives that the horizontal median $]p_i,p_{i+1}[\times\{m_j\}$ does
not meet $N(Z_{-})$. \hfill $\square$.\bigskip

{\textbf{Remark}. Similar lemmas hold with the horizontal and
vertical segments $]p_i,p_{i+1}[\times\{q_j\}$ and
$\{q_i\}\times]p_j,p_{j+1}[\,$.\medskip

\textbf{Proof of Proposition~\ref{1RZ-P0}}. The idea of the proof is
to follow the nodal set along the boundary of each square $\partial
\Sg_i$, with $i=1, 2, \ldots, r$ (say from the point $(p_i,p_i)$
anticlockwise, through $(p_{R-i},p_{i})$, $(p_{R-i},p_{R-i})$,
$(p_i,p_{R-i})$, and back to $(p_i,p_i)$), and to use the properties
of the functions $Z_{\pm}$ (no critical zeroes on the open edges, known
nodal patterns at the vertices, known interior critical zeroes, and
their localization together with Lemmas~\ref{1RZ-L2a} and
\ref{1RZ-L2b} ).

When $R=2r$ is even, it suffices to prove the result for $Z_{+}$. We
already know that the nodal set of $Z_{+}$ contains the
anti-diagonal and $\partial \Sg$. Start from $(p_1,p_1)$
horizontally. The absence of critical zero on the edge
$]0,\pi[\times\{0\}$ and Lemma~\ref{1RZ-L2a} tell us that the nodal
line can only  intertwine the edge of $\Sg_1$ untill it enters
the square $Q_{1,2r-2}$ at the point $(p_1,p_{2r-2})$. Due to the
nodal pattern in this square which contains the critical zero
$(q_1,\pi-q_1)$, the nodal line exits the square at the point
$(p_{R-1},p_2)$. By Lemma~\ref{1RZ-L2a}, and the absence of critical
zero on the edge $\{\pi\}\times]0,\pi[$, the nodal line has to
follow upwards along $x=p_{R-1}$, till the point $(p_{R-1},p_{R-1})$
where there is no choice but to get along the horizontal edge at
this point, backwards until the nodal line enters $Q_{1,R-2}$ at the
point $(p_2,p_{R-1})$. In this square, the nodal pattern in known,
and the nodal line has to leave throught the point $(p_1,p_{R-2})$
downwards along the last edge of $\Sg_1$ back to the starting point.
This is the first closed curve $\gamma_1$. We can now iterate the
procedure along $\Sg_2$, using Lemma~~\ref{1RZ-L2a} to constrain the
nodal set from both sides. After $(r-1)$ iterations, we end up with
$(r-1)$ closed curves, and we have visited every point in $\Lg$ (if
we take the diagonal into account). The curves $\gamma_i$ cannot
meet because an intersection point would be a critical zero, and we
know that the only critical zeroes of $Z_{+}$ are on the
anti-diagonal. The nodal set cannot contain any other connected
component, otherwise such a component would be entirely contained in
a white $Q$-square, and we know that this is not possible for energy
reasons. This proves Assertion (i).\medskip

We obtain the other assertions by similar arguments. \hfill
$\square$\medskip

\textbf{Remark}. We have just proved that the nodal patterns
of $Z_{\pm}$ are as suggested by the pictures (see
Figures~\ref{Plate-1-8}, \ref{Plate-1-9-1},  and \ref{Plate-1-9-2}).

\subsection{Deformation of nodal patterns}\label{SS-1RD}

In this subsection, we investigate how nodal patterns of the
family of eigenfunctions $\{\Phi^{\theta}\}$ evolve when the
parameter $\theta$ varies.

\begin{lemma}\label{1RD-L1} Assume $\theta$ is not a critical
value of the parameter,  \ie does not belong to the set $\Tg
:= \Tg_o \cup \Tg_x \cup \Tg_y$. \vspace{-3mm}
\begin{enumerate}
    \item The patterns (A) and (B) in Figure~\ref{1RQ-F1}, and the patterns
    (A), (C) and (D) in Figures~\ref{1RQ-F2e} or \ref{1RQ-F2o},
    are stable in any interval $]\theta -\epsilon ,\theta + \epsilon[ \subset
    \Tg$.
    \item Let $\theta_k \in \Tg_o$ be some critical value of $\theta$.
    When $\theta$ is close to $\theta_k$ and $\theta > \theta_k$
    (\resp $\theta < \theta_k$), the pattern (C) in
    Figure~\ref{1RQ-F1} changes to one of the patterns (A) or (B)
    (\resp (B) or (A)).
    \item Let $\theta_k \in \Tg_x \cup \Tg_y$ be some critical value of $\theta$.
    When $\theta$ is close to $\theta_k$ and $\theta > \theta_k$
    (\resp $\theta < \theta_k$), the patterns (B) in
    Figures~\ref{1RQ-F2e} or \ref{1RQ-F2o} change to one of the patterns (C) or (D)
    (\resp (D) or (C)).
\end{enumerate}
\end{lemma}

\textbf{Remark}. The proof provides more information than the
above statement.\medskip

\textbf{Proof.} The proofs are similar for $R$ even and $R$ odd. We
only sketch the proofs for $R$ even. \emph{Assertion (i).} Assume we
are in a square $Q_{i,j}$ which does not touch the boundary of the
square. In order to prove the first assertion, we consider the
segment $]0,\pi[ \times \{q_j\} \cap Q_{i,j}$ for the patterns
Figure~\ref{1RQ-F1}~(A), and the segment $\{q_i\} \times
]0,\pi[ \, \cap Q_{i,j}$ for the pattern Figure~\ref{1RQ-F1}~(B).
The argument is the same in the two cases. Let us consider the last
one. The function $y \to \Phi(q_i,y,\theta)$ has precisely two
simple zeroes in the interval $] p_j, p_{j+1}[\,$. The function $y \to
\Phi(q_i,y,\theta')$ will still have two simple zeroes for $\theta'$
close to $\theta$. As a matter of fact, when $\theta'$ varies, the
two arcs of nodal set become closer and closer, and eventually
touch, which occurs precisely when $\theta'$ reaches a critical
value
in $\Tg_o\,$. \\[4pt]
\emph{Assertion (ii).} We consider some critical value $\theta_k$,
and use the same segments as in the proof of the first assertion.
There are two cases for $\Phi(x,y,\theta_k)$ in the square
$Q_{i,j}$: it is non negative on the vertical segment and non
positive on the horizontal one, or vice and versa. Both cases are
dealt with in the same manner. For symmetry reasons, we can also
assume that $0 < \theta <\pi/4$, so that the nodal set meets
$Q_{i,j}$ if and only if $(-1)^{i+j} = -1\,$. For $\theta$ close to
and different from $\theta_k\,$, we write,
\begin{equation}\label{1RD-1}
\begin{split}
\Phi(q_i,y,\theta) =& \Phi(q_i,y,\theta_k)\\
& + (-1)^j \Big( (\cos\theta - \cos\theta_k) \sin q_i\,(-1)^j
\sin(Ry)\\
& - (\sin\theta - \sin\theta_k) (-1)^i \sin(Rq_i) \sin y\Big)\,.
\end{split}
\end{equation}

Assuming that $\Phi(q_i,y,\theta_k) \ge 0$ inside the square
$Q_{i,j}$ and looking at signs, we then see that $\Phi(q_i,y,\theta)
> 0$ inside $Q_{i,j}$ if either $j$ is even and $\theta < \theta_k$,
or $j$ is odd and $\theta > \theta_k$. This means that the nodal
pattern Figure~\ref{1RQ-F1}~(C) evolves to the nodal pattern
Figure~\ref{1RQ-F1}~(A) in these cases. Similarly, we write
\begin{equation}\label{1RD-2}
\begin{split}
\Phi(x,q_j,\theta) =& \Phi(x,q_j,\theta_k) \\
& + (-1)^j \Big( (\cos\theta - \cos\theta_k) \sin x\, (-1)^j
\sin(Rq_j) \\
& - (\sin\theta - \sin\theta_k) (-1)^i \sin(Rx)\, \sin q_j \Big)\,.
\end{split}
\end{equation}

Assuming that $\Phi(x,q_j,\theta_k) \le 0$ inside the square
$Q_{i,j}\,$, and looking at signs, we see that $\Phi(x,q_j,\theta) <
0$ inside $Q_{i,j}$ if either $j$ is even and $\theta > \theta_k\,$,
or $j$ is odd and $\theta < \theta_k\,$. This means that the nodal
pattern Figure~\ref{1RQ-F1} (C) evolves to the nodal pattern
Figure~\ref{1RQ-F1} (B) in these cases.

\emph{Assertion (iii).} The proof is similar to the proof of
Assertion (ii). \hfill $\square$

\subsection{Desingularization of $Z_{+}$}\label{SS-1RDZp}

In this  subsection, we study how the nodal set of the eigenfunction
$\{\Phi^{\theta}\}$ changes when $\theta$ varies in a small
neighborhood of $\pi/4$, while $(x,y)$ lies in the neighborhood of a
critical zero of the eigenfunction $Z_{+}$. Since $Z_{+}$ has no
critical zero when $R = 2r+1$, we only consider the case $R=2r$.
\medskip

Recall the results and notations of Subsection~\ref{SS-1RChe}. By
Lemma~\ref{1RZ-L1}, the critical zeros of $Z_{+}$ are the points
$(q_i,\pi-q_i), 1\le i\le R-2$. Take into account the fact that $R$
is even, and hence that $U_{R-1}$ is odd. For $(x,y)$ in the square
$Q(i) := Q_{i,2r-1-i}$ which contains the point $(q_i,\pi-q_i)\,$,
write
\begin{equation}\label{1RDZp-1}
\begin{split}
\sqrt{2}\, Z_{+}(q_i,y) & = \sin q_i\,\sin  y\, (-1)^{i}|U_{R-1}(\cos y) + M_i|\,,\\
\sqrt{2}\, Z_{+}(x,\pi-q_i) & = \sin q_i \,\sin x  (-1)^{i+1}|U_{R-1}(\cos x) - M_i|\,,\\
\end{split}
\end{equation}
where $M_i$ is defined in Properties~\ref{1RChe-P1}(v). These
equations give the local nodal pattern for the eigenfunction $Z_{+}$
in the square $Q(i)$. When $i$ is odd, \resp even, the nodal pattern
is given by Figure~\ref{1RDZp-F1}~(i), \resp by
Figure~\ref{1RDZp-F1}~(ii).

\begin{figure}[!ht]
\begin{center}
\includegraphics[width=10cm]{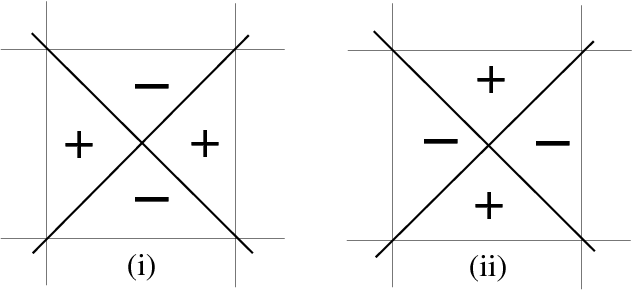}
\caption{Local nodal patterns at an interior critical
zero}\label{1RDZp-F1}
\end{center}
\end{figure}

\begin{figure}[!ht]
\begin{center}
\includegraphics[width=10cm]{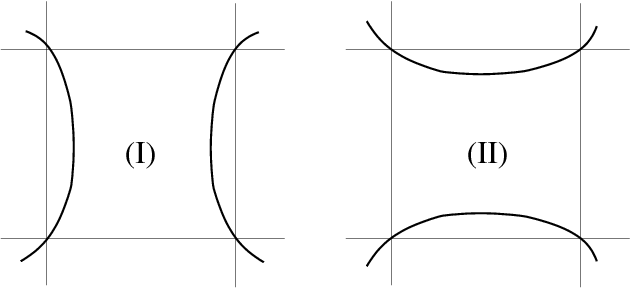}
\caption{Local nodal patterns in the absence of critical
zero}\label{1RDZp-F2}
\end{center}
\end{figure}

On the other hand, we can write,
\begin{equation}\label{1RDZp-2y}
\begin{split}
\Phi(q_i,y,\theta) & = \sin q_i\,\sin y\, \left\{ \cos\theta\,
\sqrt{2}\, Z_{+}(q_i,y) + \big( \sin\theta - \cos\theta \big)\, M_i\right\} \\
 & = (-1)^{i}\left\{ \cos\theta \, \sqrt{2}\, |Z_{+}(q_i,y)| -
 \sin q_i \,\sin y \, \big( \sin\theta - \cos\theta \big)\, |M_i|\right\}.\\
\end{split}
\end{equation}

The last factor in the second line of \eqref{1RDZp-2y} is positive
when $0 < \pi/4 - \theta \ll 1\,$. It follows that the local pattern
of the nodal set $N(\Phi)$ inside $Q(i)$, is given by
Figure~\ref{1RDZp-F2}~(II).\medskip

Similarly, we can write
\begin{equation}\label{1RDZp-2x}
\begin{split}
\Phi(x,\pi-q_i,\theta) & = \sin q_i \,\sin x \left\{ \sin\theta\,
\sqrt{2}\, Z_{+}(x,\pi-q_i) + \big( \sin\theta - \cos\theta \big)\, M_i\right\} \\
 & = (-1)^{i+1}\left\{ \sin\theta \,\sqrt{2}\, |Z_{+}(x,\pi-q_i)| +
 \sin q_i\,\sin y \, \big( \sin\theta - \cos\theta \big)\, |M_i|\right\}.\\
\end{split}
\end{equation}

The last factor in the second line of \eqref{1RDZp-2x} is positive
when $0 < \theta - \pi/4 \ll 1$. It follows that the local nodal
pattern of the nodal set $N(\Phi^{\theta})$ inside $Q(i)$, is given
by Figure~\ref{1RDZp-F2}~(I).\medskip

\textbf{Remark}. Notice that the pattern is independent of
$i\,$. When $\theta$ leaves the value $\pi/4$, all the critical zeroes
of the eigenfunction $Z_{+}$ disappear at once, and the local nodal
patterns of $\Phi^{\theta}$ in the $Q$-squares containing the
critical zeroes of $Z_{+}$ look alike, opening ``horizontally'' as
in Figure~\ref{1RDZp-F2}~(II), when $\theta < \pi/4\,$; \resp opening
``vertically'' as in Figure~\ref{1RDZp-F2}~(I), when $\theta
> \pi/4\,$, as stated in Theorem~\ref{stern-T}~(ii).

\subsection{Desingularization of $Z_{-}$}\label{SS-1RDZm}

Since $Z_{+}$ and $Z_{-}$ are symmetrical with respect to
$\{x=\pi/2\}$ when $R=2r$, we only have to consider the case
$R=2r+1$. When $R=2r+1$, the interior critical zeroes of $Z_{-}$ are
the point $(\pi/2,\pi/2)$ and the points $(q_i,q_i)$,
$(q_i,\pi-q_i)$, for $1\le i \le R-2$. Due to the symmetries with
respect to $\{x=\pi/2\}$ and $\{y=\pi/2\}$, it suffices to consider
the points $(q_i,q_i)\,$, $1 \le i \le r\,$, and the square
$Q_{i,i}$.\medskip

We can write
\begin{equation}\label{1RDZm-1}
\begin{split}
\sqrt{2}\, Z_{-}(q_i,y) & = \sin q_i \,\sin y\,  (-1)^{i+1}\; |U_{R-1}(\cos y) - M_i|\,,\\
\sqrt{2}\, Z_{-}(x,q_i) & = \sin q_i\,\sin x \,  (-1)^{i}\; |U_{R-1}(\cos x) - M_i|\,.\\
\end{split}
\end{equation}
These equations give the local nodal pattern for the eigenfunction
$Z_{-}$ in the square $Q_{i,i}$. When $i$ is even, \resp odd, the
pattern is given by Figure~\ref{1RDZp-F1}~(i), \resp by
Figure~\ref{1RDZp-F1}~(ii). \medskip

On the other hand, we can write,
\begin{equation}\label{1RDZm-2y}
\begin{split}
\Phi(q_i,y,\theta) & = \cos\theta\,\sqrt{2}\, Z_{-}(q_i,y) +
\sin q_i \,\sin y \,\big( \sin\theta + \cos\theta \big)\, M_i \\
 & = (-1)^{i}\left\{ |\cos\theta \,\sqrt{2}\, Z_{+}(q_i,y)| +
 \sin q_i \,\sin y \, \big( \sin\theta + \cos\theta \big)\, |M_i|\right\}\,.\\
\end{split}
\end{equation}

The last factor in the second line of \eqref{1RDZm-2y} is positive
when $0 < 3\pi/4 - \theta \ll 1$. It follows that the local pattern
of the nodal set $N(\Phi^{\theta})$, is given by
Figure~\ref{1RDZp-F2}~(I).\medskip

Similarly, we can write
\begin{equation}\label{1RDZm-2x}
\begin{split}
\Phi(x,q_i,\theta) & = -\sin\theta\, \sqrt{2}\, Z_{-}(x,q_i) +
\sin q_i \,\sin x \, \big( \sin\theta + \cos\theta \big)\, M_i \\
 & = (-1)^{i+1}\left\{ \sin\theta \, \sqrt{2}\, |Z_{-}(q_i,y)| -
 \sin q_i \,\sin y \, \big( \sin\theta + \cos\theta \big)\, |M_i|\right\}\,.\\
\end{split}
\end{equation}

The last factor in the second line of \eqref{1RDZm-2x} is positive
when $0 < \theta - 3\pi/4 \ll 1\,$. It follows that the local pattern
of the nodal set $N(\Phi^{\theta})$, is given by
Figure~\ref{1RDZp-F2}~(II).\medskip

We point out that the pattern is independent of the sign of $i$.
When $\theta$ leaves the value $3\pi/4$, all the critical zeroes of
the eigenfunction $Z_{-}$ disappear at once, and the local nodal
patterns of $\Phi^{\theta}$ in the squares containing the critical
zeroes of $Z_{-}$ look alike, opening ``vertically'' as in
Figure~\ref{1RDZp-F2}~(I), when $\theta < 3\pi/4\,$; \resp opening
``horizontally'' as in Figure~\ref{1RDZp-F2}~(II),  when $\theta
> 3\pi/4\,$.

\subsection{The nodal pattern of $\Phi^{\theta}$ for $\theta$
close to $\pi/4$ and $R$ even}\label{SS-1RDe}

By Properties~\ref{1RSym-P1}~(iv), we can assume that $\theta \in
[0,\pi/4]\,$. We know from Properties~\ref{1RCZV-P1}, \ref{1RCZE-P1}
and \ref{1RCZI-P1} that for $R = 2r$ and $0 < \pi/4 - \theta \ll 1\,$,
the eigenfunction $\Phi^{\theta}$ has no interior critical zero, two
non degenerate edge critical zeroes, respectively in the intervals
$\{0\}\times]p_{R-1}, \pi[$ and $\{\pi\}\times]0,p_{1}[\,$, and that
the vertices are non degenerate critical zeroes.\medskip

Using Lemma~\ref{1RD-L1}, the nodal pattern of $N(\Phi^{\theta})$
for $\theta$ close to $\pi/4\,$, is the same as the nodal pattern of
$Z_{+}$ in the $Q$-squares without critical zero, namely the white
$Q$-squares which do not meet the anti-diagonal. To determine the
nodal set of $N(\Phi^{\theta})$, it suffices to know the local nodal
patterns in the white $Q$-squares covering the anti-diagonal. This
is given by Subsection~\ref{SS-1RDZp}, \emph{the crosses at an interior
critical zero open up horizontally}, and by the description of the
critical zeroes on the vertical edges near the vertices $(0,\pi)$
and $(\pi,0)$. It is now clear that the nodal set of $\Phi^{\theta}$
for $\theta$ close enough to $\pi/4$ is connected, and divides the
square into two connected components.\medskip

\textbf{Remark}. As the diaporama of subfigures in
Figure~\ref{Plate-1-8} shows, there is another way to obtain
examples of eigenfunctions with exactly two domains, using a
deformation of one of the simplest product eigenfunctions.  The
following figures display the nodal sets for the
eigenfunctions $\Phi^{\theta}_{1,8}$ and $\Phi^{\theta}_{1,9}\,$. The
values of $\theta$ appear in the title. The values with more than
two digits correspond to the critical values of the parameter,
\ie the values of $\theta$ for which critical zeroes or
equivalently multiple points appear/disappear in the nodal set, see
\eqref{1-RChe-EX8} and \eqref{1-RChe-EX9}. The values with two
digits are intermediate values between two consecutive critical
values. The topology of the nodal set does not change in such
intervals.


\begin{figure}[!ht]
\begin{center}
\includegraphics[width=15cm]{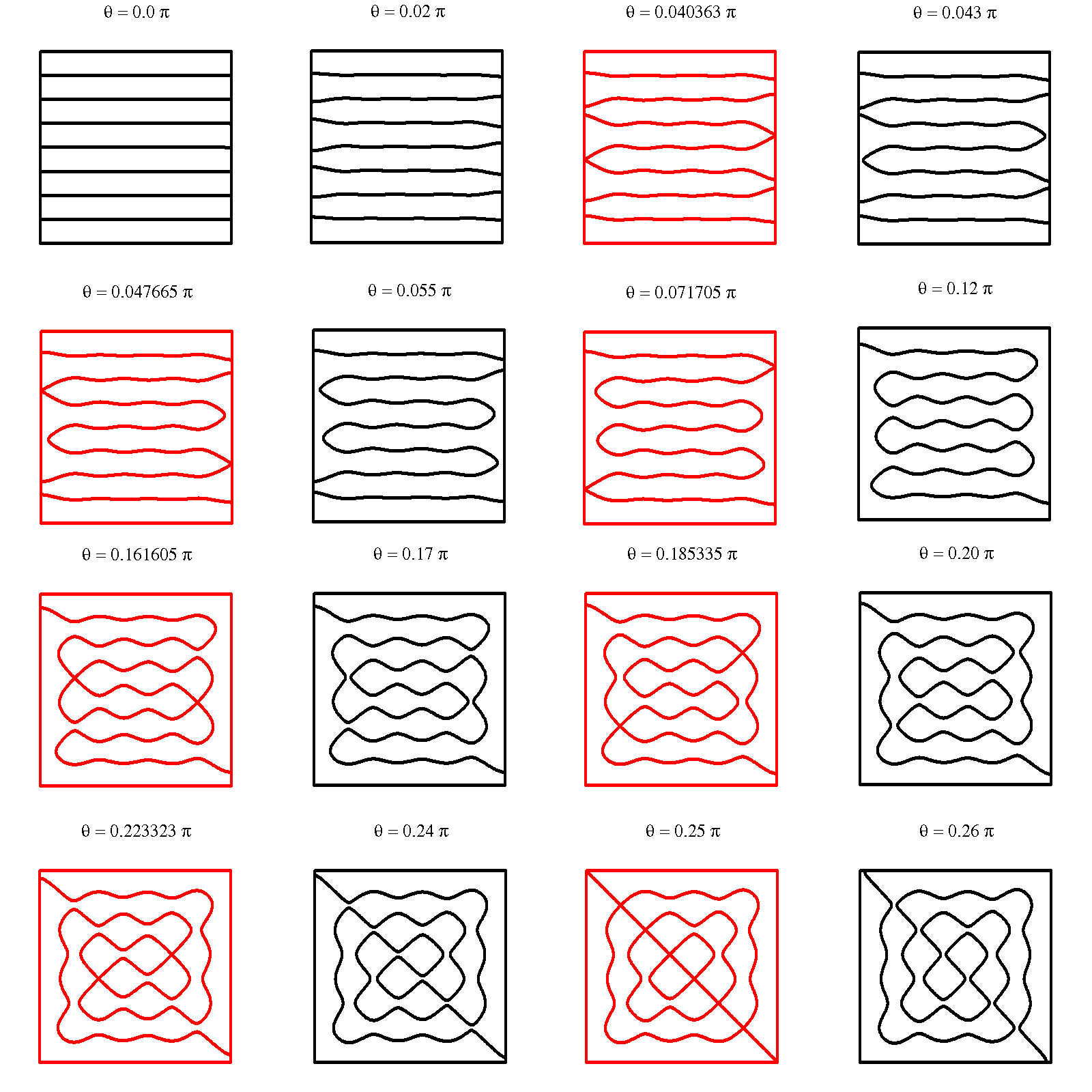}
\caption{Typical nodal patterns for the eigenvalue
$(1,8)$}\label{Plate-1-8}
\end{center}
\end{figure}


\begin{figure}[!ht]
\begin{center}
\includegraphics[width=15cm]{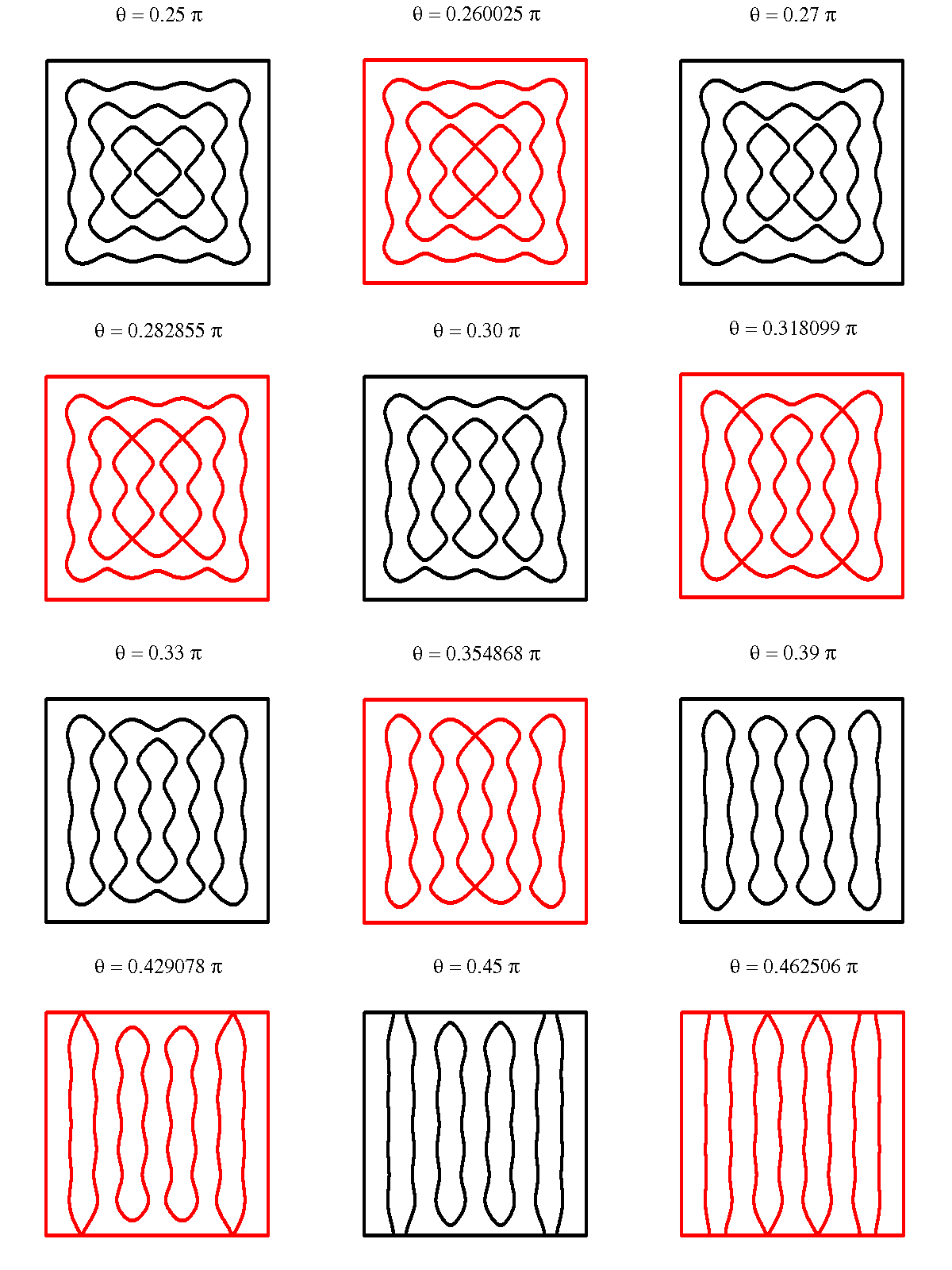}
\caption{Typical nodal patterns for the eigenvalue
$(1,9)$}\label{Plate-1-9-1}
\end{center}
\end{figure}


\begin{figure}[!ht]
\begin{center}
\includegraphics[width=15cm]{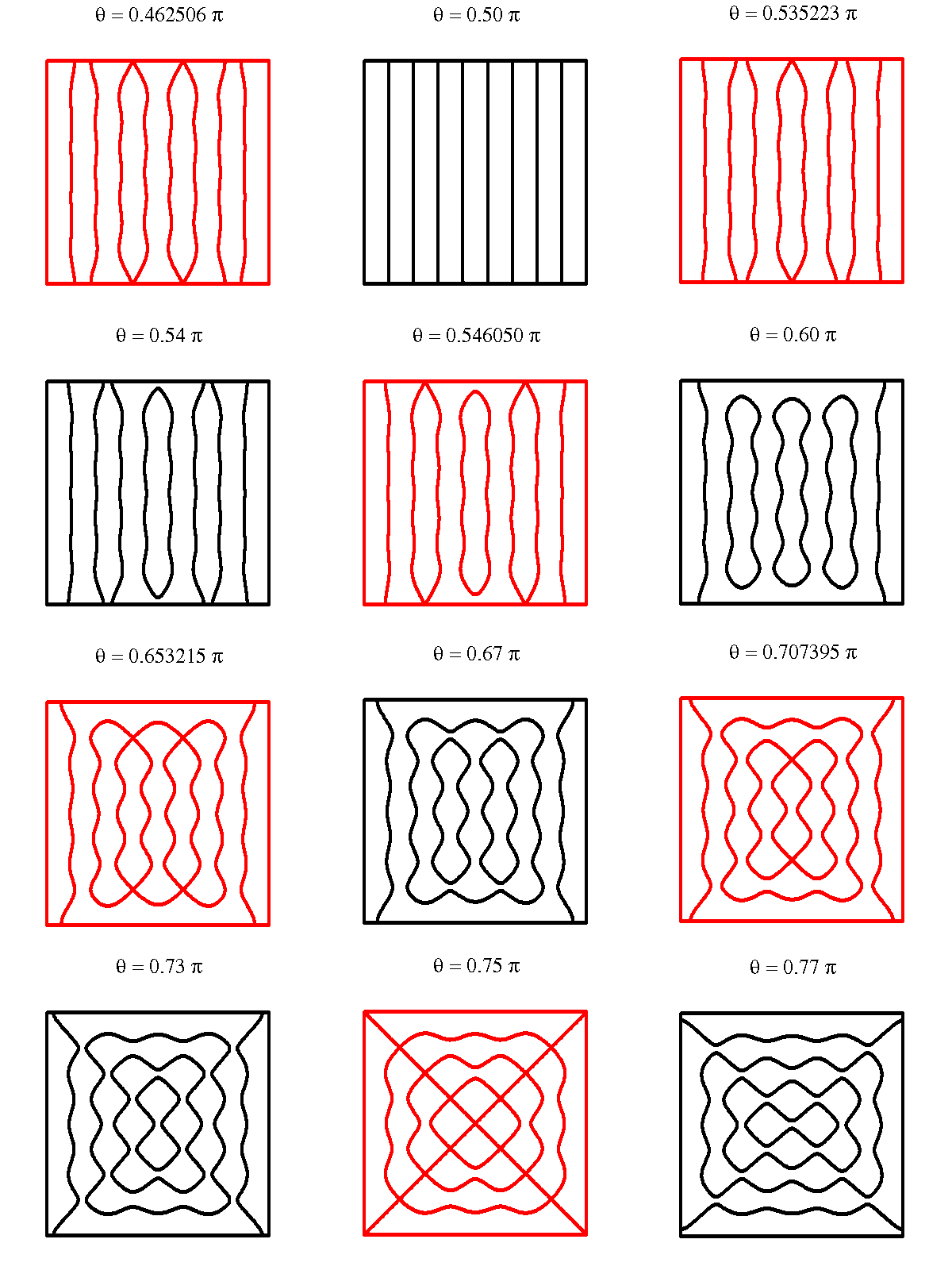}
\caption{Typical nodal patterns for the eigenvalue $(1,9)$,
continued}\label{Plate-1-9-2}
\end{center}
\end{figure}

\clearpage


\bibliographystyle{plain}

\end{document}